\documentclass[11pt]{article}
\usepackage{amsmath,amsthm,amsfonts,amssymb,amscd}

 \usepackage{amsmath}
\usepackage{enumerate} 
\usepackage{color}
\usepackage{float}
\usepackage{soul}
\usepackage{graphicx}
\usepackage[colorlinks=true]{hyperref}
\usepackage{mathpazo}
\usepackage{fancyhdr}
\usepackage{empheq}
\usepackage[margin=1in]{geometry}
\def\disp{\displaystyle}

\def\XX{\mathbb{X}}
\def\YY{\mathbb{Y}}

\def\Limsup{\mathop{{\rm Lim}\,{\rm sup}}}
\def\Liminf{\mathop{{\rm Lim}\,{\rm inf}}}

\def\tto{\rightrightarrows}
\def\Hat{\widehat}

\def\Bar{\overline}
\def\ra{\rangle}
\def\la{\langle}
\def\ve{\varepsilon}
\def\B{\mathbb{B}}
\def\h{\hfill\Box}
\def\R{\mathbb{R}}
\def\ox{\bar{x}}
\def\OX{\Bar{X}}
\def\oy{\bar{y}}
\def\OY{\Bar{Y}}

\def\ov{\bar{v}}
\def\OV{\Bar{V}}

\def\OW{\Bar{W}}

\def\OU{\Bar{U}}

\def\ri{\mbox{\rm ri}\,}

\def\Im{\mbox{\rm Im}\,}

\def\epi{\mbox{\rm epi}\,}
\def\rank{\mbox{\rm rank}\,}
\def\Tr{\mbox{\rm Tr}\,}
\def\dim{\mbox{\rm dim}\,}
\def\span{\mbox{\rm span}\,}

\def\dom{\mbox{\rm dom}\,}
\def\Ker{\mbox{\rm Ker}\,}

\def\cl*co{\mbox{\rm cl}^*\mbox{\rm co}\,}
\newcommand{\eqdef}{\stackrel{\text{\rm\tiny def}}{=}}
\def\cl{\mbox{\rm cl}\,}

\def\h{\hfill\triangle}
\def\dn{\downarrow}
\def\O{\Omega}

\def\TT{\mathbb{T}}
\def\ph{\varphi}

\def\st{\stackrel}
\def\oR{\Bar{\R}}

\def\lm{\lambda}
\def\gg{\gamma}
\def\dist{{\rm dist}\,}

\def\al{\alpha}

\def\Lm{\Lambda}

\def\MM{\mathcal{M}}

\def\hs7{\hspace*{7pt}}

\def\Id{\mathbb{I}}

\renewcommand{\theequation}{\thesection.\arabic{equation}}

\def\h{\hfill\Box}
\def\kk{\kappa}
\begin{document}

\newtheorem{Theorem}{Theorem}[section]
\newtheorem{Conjecture}[Theorem]{Conjecture}
\newtheorem{Proposition}[Theorem]{Proposition}
\newtheorem{Remark}[Theorem]{Remark}
\newtheorem{Lemma}[Theorem]{Lemma}
\newtheorem{Corollary}[Theorem]{Corollary}
\newtheorem{Definition}[Theorem]{Definition}
\newtheorem{Example}[Theorem]{Example}
\newtheorem{Fact}[Theorem]{Fact}
\newtheorem*{pf}{Proof}
\renewcommand{\theequation}{\thesection.\arabic{equation}}
\normalsize
\normalfont
\medskip
\def\endproof{$\h$\vspace*{0.1in}}
\title{\bf Geometric characterizations for strong minima with applications to nuclear norm minimization problems}
\date{}
\author{Jalal Fadili\footnote{Normandie Universit\'e, ENSICAEN, UNICAEN, CNRS, GREYC, France; email: jalal.fadili@ensicaen.fr}\,, \;\; Tran T. A. Nghia\footnote{Department of Mathematics and Statistics, Oakland University, Rochester, MI 48309, USA; email: nttran@oakland.edu}\,,\; \; and \;\; Duy Nhat Phan\footnote{Department of Mathematics and Statistics, University of Massachusetts Lowell, Lowell, MA 01854, USA; duynhat\_phan@uml.edu} }

\maketitle

\begin{quote}
{\small \noindent {\bf Abstract.} In this paper, we introduce several geometric characterizations for strong minima of optimization problems. Applying these results to nuclear norm minimization problems allows us to obtain new necessary and sufficient quantitative conditions for this important property. Our characterizations for strong minima are  weaker than the Restricted Injectivity and Nondegenerate Source Condition, which are usually used to identify solution uniqueness of nuclear norm minimization problems. Consequently, we obtain the minimum (tight) bound on the number of measurements for (strong) exact recovery of low-rank matrices.  }


\end{quote}
\maketitle

\section{Introduction}
\setcounter{equation}{0}
Strong minima is an important property at a local minimizer of an optimization problem such that the difference between the cost value and the optimal value is bigger than a proportional of the norm square of the difference between the corresponding feasible solution and the minimizer. It is an error bound condition that has various applications to sensitivity analysis, robustness, and complexity of algorithms \cite{BT80,BZ82,BCS99,BS00,CDZ17,FM68, FMP18,FNT21,LFP17,I79,R88,SW99}.  Finding necessary and sufficient second order conditions for {\em strong minima} is a classical research area. For nonlinear programming, the first results in this direction were probably established in \cite{FM68} under some restrictive conditions. Complete second order characterizations for nonlinear programming under mild conditions such as the Mangasarian-Fromovitz constraint qualification were obtained later in \cite{BT80,I79}. For constrained (nonpolyhedral) optimization problems with smooth data, necessary and sufficient second order conditions for strong minima are much more involved. They often contain nontrivial ``sigma terms", which represent curvatures of some nonpolyhedral structures in the problem. Another important feature of these conditions is that they are usually formulated  as ``minimax'' conditions  such that {\em Lagrange multipliers} are dependent on the choice of vectors in the critical cone; see, e.g.,  \cite{BCS99,BS00}. Although the aforementioned sigma terms are fully calculated for many seminal classes of optimization such as semi-infinite programming, semi-definite programming, and second order cone programming, their calculations are complicated in general. Moreover, checking the (minimax) sufficient second order conditions is quite a hard task numerically. 

An important problem that motivates our study in this paper is the {\em nuclear norm minimization problem} 
\begin{equation}\label{Nu0}
    \min_{X\in \R^{n_1\times n_2}}\quad \|X\|_*\quad \mbox{subject to}\quad \Phi X=M_0,
\end{equation}
where $\|X\|_*$ is the nuclear norm of an $n_1\times n_2$ matrix $X$, $\Phi:\R^{n_1\times n_2}\to \R^m$ is a linear operator, and $M_0$ is a known vector (observation) in $\R^m$. This problem is considered the tightest convex relaxation of the celebrated NP-hard {\em affine rank minimization problem} with various  applications in computer vision, collaborative filtering, and data science; see, e.g., \cite{ALMT14,CRPW12,CR09,CR13,RFP10}. There are  some essential  reasons to study strong minima of this problem. First, strong minima of problem \eqref{Nu0} guarantees the linear convergence of some proximal algorithms for solving problem \eqref{Nu0} and related problems; see, e.g., \cite{CDZ17,LFP17,ZS17}. Second, it is also sufficient for solution uniqueness and robustness of problem \eqref{Nu0} \cite{FMP18,FNT21}. Solution uniqueness for problem \eqref{Nu0} is a significant property in recovering the original low-rank solution $X_0\in \R^{n_1\times n_2}$ from observations $M_0=\Phi X_0$. In \cite{CR09,CR13}, Cand\`es and Recht introduced a {\em nondegenerate condition} sufficient for solution uniqueness of  problem \eqref{Nu0}.  It plays an important role in their results about finding a small bound for measurements $m$  such that solving problem \eqref{Nu0} recovers exactly $X_0$ from observations $M_0$ over a Gaussian linear operator $\Phi$. Their condition is recently revealed in \cite{FNT21} to be a complete characterization for the so-called {\em sharp minima} introduced by  Crome \cite{C77} and Polyak \cite{P79} independently. This special property of problem \eqref{Nu0} at $X_0$ guarantees the {\em robust recovery} with a linear rate in the sense that  any solutions of the following low-rank optimization problem 
\begin{equation}\label{LSN}
\frac{1}{2}\|\Phi X-M\|^2+\mu\|X\|_*
\end{equation}
converge to $X_0$ with a linear rate as $\mu\dn 0$ provided that $\|M-M_0\|\le c\mu$ for some constant $c>0$; see \cite{CP10}. When strong minima occurs in problem \eqref{Nu0}, \cite{FNT21} shows that the convergence rate is  H\"olderian with order $\frac{1}{2}$. It is also worth noting that solution uniqueness of nuclear norm minimization problem \eqref{Nu0} can be characterized geometrically via the {\em descent cone} \cite{ALMT14,CRPW12}. As the descent cone is not necessarily closed, using it to check solution uniqueness numerically is not ideal. Another geometric characterization for solution uniqueness of problem \eqref{Nu0} is established recently in  \cite{HP23}, but the set in their main condition is  not closed too. As strong minima is necessary for sharp minima and sufficient for solution uniqueness, it is open to study the impact of strong minima to exact recovery.

Sufficient second order condition for strong minima of problem \eqref{Nu0} can be obtained from \cite[Theorem~12]{CDZ17}. The approach in  \cite{CDZ17} is rewriting problem \eqref{Nu0} as a {\em composite optimization problem} and applying the classical results in \cite{BCS99,BS00}. Some second order analysis on {\em spectral functions} including the nuclear norm studied in \cite{CD19, CDZ17,D17, MS23,ZZX13} could be helpful in understanding this result, but these second order computations applied on the nuclear norm still look complicated. Most importantly, the sufficient second order condition obtained in \cite[Theorem~12]{CDZ17} is still in a minimax form, which makes it hard to check. Our main questions throughout the paper are: 1.  Is it possible to obtain simple necessary and sufficient conditions for strong minima of problem \eqref{Nu0}? 2. Can we  avoid the minimax form usually presented in these kinds of second order sufficient conditions? and 3. Is there any efficient way to check strong minima of problem \eqref{Nu0} numerically? 

{\bf Our contribution.} To highlight the new ideas and direct to the bigger picture, we study in Section~3 the following composite optimization problem 
\begin{equation}\label{Comp}
    \min_{x\in \XX}\quad f(x)+g(x),
\end{equation}
where $\XX$ is a finite dimensional space,  $f:\XX\to \R$ is a twice continuously differentiable function,  and $g:\XX\to \oR\eqdef\R\cup\{+\infty\}$ is a proper lower semi-continuous (nonsmooth) function. It covers problem \eqref{LSN} and many modern ones in optimization. One of the most popular ways to characterize strong minima of problem \eqref{Comp} is using the  {\em second subderivative} \cite{BS00,R88,RW98}. As the function $g$ is nonsmooth, second subderivative of $g$ is hard to compute in general. To avoid this computation, we assume additionally that the function $g$ satisfies the classical {\em quadratic growth condition} \cite{BI95a,W94}. In the case of convex functions, it is shown in Section~3 that $\ox$ is a strong solution of problem \eqref{Comp} if and only if $\ov\eqdef-\nabla f(\ox)\in \partial g(\ox)$ and the following {\em geometric} condition holds
\begin{equation}\label{GEO}
\Ker \nabla^2 f(\ox)\cap T_{(\partial g)^{-1}(\ov)}(\ox)=\{0\}, 
\end{equation}
where $\partial g:\XX\tto \XX^*$ is the  subdifferential mapping of $g$, $\Ker \nabla^2 f(\ox)$ is the nullspace of the Hessian matrix $\nabla^2 f(\ox)$, and $T_{(\partial g)^{-1}(\ov)}(\ox)$ is the Bouligand {\em contingent cone} at $\ox$ to  $(\partial g)^{-1}(\ov)$. The action of the contingent cone to a first order structure in the above condition tells us that \eqref{GEO} is actually a second order condition. But it looks simpler for computation than the second subderivative; see, e.g., our Corollary~\ref{CoNu} for the case of nuclear norm. Our results also work when $g$ is not convex; see Section~3 for more detailed analysis. 

Another problem considered in Section~3 is the following convex optimization problem with linear constraints
\begin{equation}\label{ConP}
    \min_{x\in \XX}\quad g(x)\quad  \mbox{subject to}\quad \Phi x\in K,
\end{equation}
where $g:\XX\to \oR$ is a continuous (nonsmooth) convex function, $\Phi:\XX\to \YY $ is a linear operator between two finite dimensional spaces, and $K$ is a closed  polyhedral set of $\YY$. This problem covers the nuclear norm minimization problem \eqref{Nu0} and a handful of other significant optimization problems \cite{ALMT14,AV11,CRPW12}. When $g$ is twice continuously differentiable, characterizations for strong minima of problem \eqref{ConP} are  simple; see, e.g., \cite[Theorem~3.120]{BS00}. But when $g$ is not differentiable, characterizing strong minima is much more involved. A standard approach is rewriting problem \eqref{ConP} as a composite problem, which can be represented by a constrained optimization problem with smooth data \cite{BCS99,BS00, I91}.  To obtain similar geometric characterization to \eqref{GEO} for problem \eqref{ConP}, we additionally assume the function $g$ satisfies both quadratic growth condition and {\em second order regularity}. The latter condition was introduced by Bonnans, Cominetti, and Shapiro in \cite{BCS99} to close the gap between necessary and sufficient second order  optimality conditions for constrained and composite optimization problems. But the extra assumption of quadratic growth condition  on the function $g$ in this paper allows us to achieve geometric characterizations for strong minima of problem \eqref{ConP} in Theorem~\ref{thm5}. Many vital classes of nonsmooth functions $g$ satisfy both the quadratic growth condition and second order regularity. To list a few, we have classes of  piecewise linear-quadratic convex functions \cite{RW98},  spectral functions \cite{CDZ17,MS23}, $\ell_1/\ell_2$ norms \cite{ZS17}, and  indicator functions to the set of  positive semi-definite matrices \cite{CST19}. Our results are applicable not only to  nuclear norm minimization problem \eqref{Nu0}, but also other different  problems in optimization. 

As the nuclear norm satisfies both the quadratic growth condition and second order regularity,  geometric characterizations for strong minima of low-rank problems \eqref{Nu0} and \eqref{LSN} are foreseen. But our studies on problems \eqref{Nu0} and \eqref{LSN} in Section 4 and 5 are not just straightforward applications. In Section~4, we derive simple calculation of the second order structure in  \eqref{GEO} for the case of nuclear norm. Furthermore, some quantitative characterizations for strong minima of problem \eqref{LSN} are obtained via the so-called Strong Restricted Injectivity and Analysis Strong Source Condition that inherit the same terminologies introduced recently in \cite{FNT21} to characterize strong minima/solution uniqueness of group-sparsity optimization problems. These conditions are weaker than the well-known Restricted Injectivity and Nondegenerate Source Condition used in \cite{CR09,CR13} as sufficient conditions for solution uniqueness of nuclear norm minimization problem \eqref{Nu0}; see also \cite{GHS11} for the case of $\ell_1$-norm. Both conditions can be verified numerically. In Section 5, we obtain new characterizations for strong minima of  problem \eqref{Nu0}. Our conditions are not in the form of minimax problems. Indeed, Theorem~\ref{SSN} shows that $X_0$ is a strong solution of problem \eqref{Nu0} if and only if there exists a dual certificate $\OY\in \Im\Phi^*\cap \partial \|X_0\|_*$ such that
\[
\Ker \Phi \cap T_{(\partial \|\cdot\|_*)^{-1}(\OY)}(X_0)=\{0\},
\]
which  has some similarity with \eqref{GEO}. Necessary and sufficient conditions for strong minima obtained in this section reveal some interesting facts about exact recovery for nuclear minimization problem \eqref{Nu0}. For example, one needs at least $\frac{1}{2}r(r+1)$ measurements for $M_0=\Phi X_0$ to recover exactly  the matrix $X_0$ of rank $r$ as a strong solution of problem \eqref{Nu0}. This bound for $m$ is very small, but it is tight in the sense that we can construct infinitely many linear operators $\Phi:\R^{n_1\times n_2}\to \R^{\frac{1}{2}r(r+1)}$  such that solving problem \eqref{Nu0} recovers exactly $X_0$. Another compelling result in Section 5 shows that  the {\em low-rank representation problem} \cite{HO14,LLYS12} always has strong minima, when the linear operator $\Phi$ in  \eqref{Nu0} is any $q\times n_1$ matrix. 

Finally in this paper, we discuss numerical methods to check strong minima and compare the results with sharp minima and  solution uniqueness. For example, with 100 nuclear norm minimization problems via standard Gaussian linear operators $\Phi$ and $460$  measurements $M_0\in \R^{460}$ observed from the original matrix $X_0\in \R^{40\times 40}$ of rank $3$, the case of exact recovery is about 80\%, in which  about 40\% problems have sharp minima and the other 40\% problems have strong minima. As the traditional approach for exact recovery in \cite{ALMT14,CRPW12,CR13} is via  sharp minima \cite{FNT21}, seeing more unique (strong, non-sharp) solutions in these numerical experiments gives us  a more complete picture about exact recovery when number of measurements $m$ is not big enough; see our Section~6 for further  numerical experiments.

\section{Preliminaries}
\setcounter{equation}{0}

Throughout this paper, we suppose that $\XX$ is an Euclidean space with norm $\|\cdot\|$ and  $\XX^*$ is its dual space endowed with the inner product $\la v,x\ra$ for any $v\in \XX^*$ and $x\in \XX$. $\B_r(\ox)$ is denoted by the closed ball with center $\ox\in \XX$ and radius $r>0$. Let  $\ph:\XX\to \oR\eqdef\R\cup\{+\infty\}$ be a proper extended real-valued function with nonempty domain $\dom \ph\eqdef\{x\in \XX|\; \ph(x)<\infty\}\neq \emptyset$. A point $\ox\in \dom \ph$ is called a \emph{strong solution} (or strong minimizer) of $\ph$ if there exist $c,\ve>0$ such that 
\begin{equation}\label{Stro}
    \ph(x)-\ph(\ox)\ge c\|x-\ox\|^2\quad \mbox{for all}\quad x\in \B_\ve(\ox). 
\end{equation}  
In this case, we say that strong minima occurs at $\ox$. 
The study of strong minima is usually based on second order theory, where the following structures of subderivatives play crucial roles; see, e.g., \cite{BZ82,BS00,I91,R88,RW98}. 
\begin{Definition}[Subderivatives] For a function $\ph:\XX\to \oR$ and $\ox\in \dom \ph$,  the {\em subderivative} of $\ph$ at $\ox$ is the function  $d\ph(\ox):\XX\to \oR$ defined by
\begin{equation}\label{subd1}
d\ph(\ox)(w)\eqdef\liminf_{t\dn 0, w^\prime\to w}\dfrac{\ph(\ox+tw^\prime)-\ph(\ox)}{t}\quad \mbox{for} \quad w\in X.
\end{equation}
The {\em second subderivative} of $\ph$ at $\ox$ for $\ov\in \XX^*$ is the function $d^2 \ph(\ox|\ov):\XX\to \oR$  defined by 
\begin{equation}\label{Subd}
d^2\ph(\ox|\ov)(w)\eqdef\liminf_{t\dn 0, w^\prime \to w}\dfrac{\ph(\ox+tw^\prime)-\ph(\ox)-t\la \ov, w^\prime\ra}{\frac{1}{2}t^2}\quad \mbox{for} \quad w\in X.
\end{equation}
The {\em parabolic subderivative} of $\ph$ at $\ox$ for $w\in \dom d\ph(\ox)(\cdot)$ with respect to $z\in \XX$ is defined by
\begin{equation}\label{SED}
d^2\ph(\ox)(w|z)\eqdef\liminf_{t\dn 0, z^\prime \to z}\dfrac{\ph(\ox+tw+\frac{1}{2}t^2z^\prime)-\ph(\ox)-td\ph(\ox)(w)}{\frac{1}{2}t^2}.
\end{equation}
\end{Definition}

Parabolic subderivatives were introduced by Ben-Tal and Zowe  in \cite{BZ82}  to study strong minima; see also \cite[Theorem~13.66]{RW98}. Second subderivatives dated back to the seminal work of Rockafellar \cite{R88} with  good calculus \cite{I91} for many important classes of functions.  It is well known \cite[Theorem~13.24]{RW98} that $\ox$ is a strong solution of $\ph$ if and only if $0\in \partial \ph(\ox)$ and 
\begin{equation}\label{SS}
    d^2\ph(\ox|0)(w)>0\quad \mbox{for all}\quad w\neq 0. 
\end{equation}
Here $\partial\ph(\ox)$ stands for the  Mordukhovich {\em limiting subdifferential} of $\ph$ at $\ox$ \cite{M1}:
\begin{equation}\label{LimS}
    \partial \ph(\ox)=\left\{
v\in \XX^*|\; \exists (x_k,v_k)\st{X\times \XX^*}\longrightarrow (\ox,v),\, \liminf_{x\to x_k}\dfrac{\ph(x)-\ph(x_k)-\la v_k,x-x_k\ra}{\|x-x_k\|}\ge 0\right\}.
\end{equation}
When $\ph$ is a proper l.s.c. convex function, this subdifferential coincides with the subdifferential in the classical convex analysis
\begin{equation}\label{Cv}
  \partial \ph(\ox)=\left\{
v\in \XX^*|\; \ph(x)-\ph(\ox)\ge \la v,x-\ox\ra, x\in \XX\right\}.  
\end{equation}
We denote $\ph^*:\XX^*\to \oR$ by the \emph{Fenchel conjugate} of $\ph$:
\begin{equation}\label{Fech}
\ph^*(v)\eqdef\sup\{\la v,x\ra-\ph(x)\ra|\;x\in \XX\}\quad \mbox{for}\quad v\in \XX^*. 
\end{equation}


Next let us recall here some first and second order tangent structures \cite{BS00,RW98} on a nonempty closed set $K$ of $\XX$ that widely used in this paper.



\begin{Definition}[tangent cones] Let $K$ be a closed set of $\XX$. The Bouligand \emph{contingent cone} at the point $\ox\in K$ to $K$ is defined by
\begin{equation}\label{CC}
    T_{K}(\ox)\eqdef\disp\Limsup_{t\dn 0}\frac{K-\ox}{t}=\left\{w\in \XX|\;\exists\,t_k\dn 0, w_k\to w,  \ox+t_kw_k\in K\right\}.
\end{equation}
The {\em inner and
outer second order tangent set} to $K$ at $\ox\in K$ in the direction $w\in \XX$ are defined, respectively,  by
\begin{equation}\label{Ti3}
\begin{aligned}
T^{i,2}_K(\ox|w)& \eqdef\disp\Liminf_{t\dn 0}\frac{K-\ox-tw}{\frac{1}{2}t^2}=\left\{z\in \XX|\; \forall\, t_k\dn 0, \exists\, z_k\to z, \ox+t_kw+\frac{1}{2}t_k^2z_k\in K\right\}
\quad \mbox{and} \\
T^{2}_K(\ox|w)& \eqdef\disp\Limsup_{t\dn 0}\frac{K-\ox-tw}{\frac{1}{2}t^2}=\left\{z\in \XX|\; \exists\,t_k\dn 0, z_k\to z, \ox+t_kw+\frac{1}{2}t_k^2z_k\in K\right\}.
\end{aligned}
\end{equation}

\end{Definition}
The contingent cone $T_K(\ox)$ is a closed set. It contains all $w\in \XX$ such that there exists a sequence $\{t_k\}\dn 0$ such that $\dist(\ox+t_kw;K)=o(t_k)$, where $\dist(x;K)$ denotes the distance from $x\in \XX$ to $K$:
\begin{equation}\label{Dist}
    \dist(x;K)=\min\{\|x-u\||\; u\in K\}. 
\end{equation}
Similarly, the inner second order tangent set to $K$ at $\oy$ is
\begin{equation}\label{ConT}
T^{i,2}_K(\ox|w)=\left\{z\in \XX|\; \dist(x+tw+\frac{1}{2}t^2z;K)=o(t^2), t\ge 0\right\}. 
\end{equation}
When $K$ is convex, it is well-known that 
\[
T_K(\ox)=\{w\in \XX|\; \dist(\ox+tw;K)=o(t), t\ge 0\}. 
\]
Since the function $\dist(\cdot;K)$ is a convex function, $T_K(\ox)$ is a convex set. In this case, the inner second order tangent set $T^{i,2}_K(\oy|w)$ is also convex due to the same reason and formula \eqref{ConT}.  Moreover, the dual of contingent cone is the \emph{normal cone} to $K$ at $\ox$:
\begin{equation}\label{Normal}
    N_K(\ox)\eqdef [T_K(\ox)]^-=\{v\in \XX^*|\; \la v, x-\ox\ra\le 0\quad \mbox{for all}\quad x\in K\}.
\end{equation}
It  is also the subdifferential of the indicator function $\iota_K$ to the set $K$, which is defined by $\iota_K(x)=0$ if $x\in K$ and $+\infty$ otherwise. The normal cone can be characterized via the \emph{support function} to $K$
\begin{equation}\label{Supp}
    \sigma_K(v)\eqdef\sup\{\la v,x\ra|\; x\in K\}\quad \mbox{for all}\quad v\in \XX^*
\end{equation}
with $N_K(\ox)=\{v\in \XX^*|\; \sigma_K(v)\le \la v, \ox\ra\}$. 

To characterize  strong minima for constrained optimization problems,  Bonnans, Cominetti, and Shapiro \cite[Definition~3]{BCS99} introduced the following  {\em second order regular condition} on $K$; see also \cite[Definition~3.85]{BS00}.

\begin{Definition}[Second order regularity]\label{RS} The set $K$ is called to be second order regular at $\ox\in K$ if for any $w\in T_K(\ox)$ the outer second order tangent set  $T^2_K(\ox|w)$ coincides with the inner second order tangent set $T^{i,2}_K(\ox|w)$ and for any sequence $x_k\in K$ of the form $x_k=\oy+t_kw+\frac{1}{2}t_k^2r_k$
\[
\lim_{k\to \infty}\dist(r_k;T^2_K(\ox,w))=0.
\]
The proper l.s.c. function $\ph:\XX\to \oR$ is said to be second order regular at $\ox\in \dom \ph$ if its epigraph of $\ph$ is  second order regular at $(\ox,\ph(\ox))$.
\end{Definition}



The class of second order regular sets cover many important sets in optimization such as any polyhedral set and the set of positive semi-definite matrices and the second order ice cream cone; see, e.g.,  \cite{BS00}. Piecewise linear quadratic convex functions are  second order regular \cite{BCS99}.   Recently, it is proved in \cite{CDZ17} some special spectral functions are also second order regular.

When the function $\ph:\XX\to \oR$ is l.s.c. convex and second order regular at $\ox\in \dom \ph$, we note from \cite[Proposition~3.41]{BS00}
\begin{equation}\label{Tepi}
    T^{i,2}_{{\rm epi}\, \ph}((\ox,\ph(\ox))|(w,d\ph(\ox)(w)))=\epi d^2\ph(\ox)(w| \cdot)
\end{equation}
for any $w\in \dom d\ph(\ox)$. This is a convex set, which implies that $d^2\ph(\ox)(w| \cdot)$ is a convex function.  In this case, it is known from \cite[Proposition~103]{BS00} that $\ph$ is {\em parabolically regular} at $\ox$ in a direction $w\in \XX$ for $v\in \XX^*$ in the sense that  
\begin{equation}\label{Bolic}
  d^2\ph(\ox|v)(w)=- [d^2\ph(\ox)(w|\cdot)]^*(v),
\end{equation}
which is the Fenchel conjugate of the function $d^2\ph(\ox)(w|\cdot)$ at $v$ provided that  the pair $(w,v)\in \XX\times \XX^*$ satisfies the condition $\la v,w\ra=d\ph (\ox)(w)$.

Next let us slightly modify \cite[Theorem~3.108 and Theorem 3.109]{BS00}, which give necessary and sufficient conditions for strong solutions of the following composite problem


%


\begin{equation}\label{CoP}
    \min_{x\in \XX}\quad g(F(x)),
\end{equation}
where $F:\XX\to \YY$ is a twice continuously differentiable mapping and $g:\YY\to \oR$ is a l.s.c. proper convex function. Suppose that  $y_0= F(x_0)\in \dom g$ with $x_0\in \XX$. The Robinson's constraint qualification at $x_0$ for this composite problem is known as 
\begin{equation}\label{Rob}
    0\in {\rm int}\, (y_0+\nabla F(x_0)\XX-\dom g);
\end{equation}
see, e.g., \cite{BCS99,BS00}. The feasible point $x_0$ is a called a {\em stationary point} of problem \eqref{CoP} if there exists a {\em Lagrange multiplier} $\lm\in \YY^*$ such that 
\begin{equation}\label{Lagr}
    \nabla F(\ox)^*\lm=0\quad \mbox{and}\quad \lm \in \partial g(y_0).
\end{equation}  
\begin{Theorem}[Second order characterizations for strong solutions of composite problems]\label{ComP} Suppose that Robinson's constraint qualification \eqref{Rob} holds at a stationary point $x_0$ and that the function $g$ is second order regular at $y_0$. Then $x_0$ is a strong solution of problem \eqref{CoP} if and only if for any nonzero  $w$ in the critical cone
\begin{equation}\label{CriC}
C(x_0)\eqdef\{u\in \XX|\; dg(y_0)(\nabla F(x_0)u)= 0\},
\end{equation}
there exists a Lagrange multiplier $\lm$ satisfying condition \eqref{Lagr} such that 
\begin{equation}\label{SOC}
    \la\lm, \nabla^2 F(x_0)(w,w)\ra+d^2g(y_0|\lm)(\nabla F(x_0)w)>0.
\end{equation}
\end{Theorem}
\noindent{\bf Proof.} Let us justify the sufficient part first. According to \cite[Theorem~3.109]{BS00}, $x_0$ is a strong solution of problem \eqref{CoP} provided that  for any $w\in C(x_0)\setminus\{0\}$ there exists a Lagrange multiplier $\lm \in \YY^*$ satisfying \eqref{Lagr} such that 
\begin{equation}\label{SOC1}
    \la\lm, \nabla^2 F(x_0)(w,w)\ra-\Psi^*(\lm)>0,
\end{equation}
where $\Psi(\cdot)\eqdef d^2g(y_0)(\nabla F(x_0)w|\cdot)$. Since $g$ is convex and second order regular at $y_0$, equation \eqref{Tepi} for the function $g$ tells us that $d^2g(y_0)(\nabla F(x_0)w|\cdot)$ is a convex function for any $w\in C(x_0)$. Moreover, note that 
\[
\la \lm, \nabla F(\ox)w\ra=\la \nabla F(\ox)^*\lm,w\ra=0=dg(\nabla F(\ox)w).
\]
We obtain from \eqref{Bolic} that $d^2g(y_0|\lm)(\nabla F(\ox)w)=-\Psi^*(\lm)$. This ensures the equivalence between \eqref{SOC1} and \eqref{SOC}. Thus $x_0$ is a strong solution of problem \eqref{CoP} provided that condition \eqref{SOC} holds. 

To prove the necessary part, we note again that the function $d^2g(y_0)(\nabla F(x_0)w|\cdot)$ is convex, then there is no gap between second order necessary and sufficient condition, i.e., condition \eqref{SOC1} is also a necessary condition for strong solution $x_0$; see  \cite[Theorem~5.2]{BCS99}   
 or \cite[Theorem~3.108]{BS00}. Due to the equivalence of \eqref{SOC} and \eqref{SOC1} above, condition \eqref{SOC} is also necesary for the strong minima at $x_0$.   \endproof

As described in \eqref{SOC}, the existence of Lagrange multiplier is dependent on the choice of each vector in the critical cone. Under the Robinson's constraint qualification \eqref{Rob},
\eqref{SOC} is a minimax condition in the sense that it is equivalent to 
\begin{equation}\label{MinMax}
    \min_{w\in C(x_0), \|w\|=1}\max_{\lm \in \Lm(x_0)}\big[\la \lm,\nabla^2F(x_0)(w,w)\ra+d^2g(y_0|\lm)(\nabla F(x_0)w)\big]>0,
\end{equation}
where $\Lm(x_0)$ is the set of all Lagrange multipliers satisfying \eqref{Lagr}. It is hard to check this condition numerically. On the other hand, its maximin version 
is more desirable, as it means that there may exist a Lagrange multiplier $\lm$ such that inequality \eqref{SOC} is valid for any $w\in C(x_0)\setminus\{0\}$. However, it is not clear how to close the gap between the minimax and maximin. For the case of \eqref{Nu0}, we will obtain some kind of maximin condition for strong minima in Theorem~\ref{SSN}.

\section{Geometric characterizations for strong minima of optimization problems}
\setcounter{equation}{0}
\subsection{Geometric characterizations for strong minima of unconstrained optimization problems}
In this subsection, we consider  the following composite optimization problem 
\begin{equation}\label{P}
    \disp\min_{x\in \XX}\quad \ph(x)\eqdef f(x)+g(x),
\end{equation}
where $f, g:\XX \to \oR$ are proper functions such that  ${\rm int}\, (\dom f)\cap \dom g\neq\emptyset$,  $f$ is twice continuously differentiable in ${\rm int}\, (\dom f)$, and $g$ is lower semi-continuous. We assume that $\ox\in {\rm int}\, (\dom f)\cap \dom g$ is a {\em stationary point} of problem \eqref{P} in the sense that 
\[
0\in \partial \ph(\ox)=\nabla f(\ox)+\partial g(\ox)
\]
due to the sum rule for limiting subdifferential; see, e.g.,  \cite[Proposition~1.107]{M1} or 
\cite[Exercise 10.10]{RW98}. Obviously, $\ox$ is a stationary point if and only if  $-\nabla f(\ox)\in \partial g(\ox)$.

To characterize strong minima at the stationary point $\ox$, one of the most  typical methods  is using  the second subderivative $d^2\ph(\ox|0)$ defined in  \eqref{SS}. As the function $f$ is twice continuously differentiable at $\ox$, it is well-known \cite[Exercise 13.18]{RW98} that 
\begin{equation}\label{SR2}
    d^2\ph(\ox|0)(w)=\la\nabla^2 f(\ox)w,w\ra+d^2 g(\ox|-\nabla f(\ox))(w)\quad \mbox{for}\quad w\in \XX.  
\end{equation}
Since $g$ is possibly nonsmooth in many structured optimization problems, the computation of $d^2 g(\ox|-\nabla f(\ox))(w)$ could be quite challenging. In this section, we establish several new necessary and sufficient conditions for strong minima without computing second subderivatives under an additional assumption that  the function $g$ satisfies the following \emph{quadratic growth condition} \cite{BI95a,W94}; see also \cite[Section~3.5]{BS00}.

\begin{Definition}[Quadratic growth conditions] Let $g:\XX\to \oR$ be a proper l.s.c. function and $S$ be a closed subset of $\XX$ with $\ox\in \dom g\cap S$. 
We say that $g$ satisfies the quadratic growth condition at $\ox$ for some $\ov\in \partial g(\ox)$ with respect to $S$ if there exist constants $\ve,\delta>0$ and modulus $\kk>0$ such that
\begin{equation}\label{LQG}
    g(x)-g(\ox)-\la \ov,x-\ox\ra\ge \frac{\kk}{2}[\dist(x;S)]^2\qquad\mbox{for all}\qquad  x\in \B_\ve^\delta(\ox|\ov)
\end{equation}
with  $\B_\ve^\delta(\ox|\ov)\eqdef\{x\in \B_\ve(\ox)|\; g(x)-g(\ox)-\la \ov,x-\ox\ra<\delta\}$.
The function $g$ is said to satisfy the quadratic growth condition at $\ox$ for $\ov$, if it satisfies this condition at $\ox$ for $\ov\in \partial g(\ox)$ with respect to 
\begin{equation}\label{Sol}
    S(\ox,\ov)\eqdef\left\{x\in \XX|\; g(x)-\la \ov,x\ra\le g(\ox)-\la \ov,\ox\ra\right\}. 
\end{equation}
Finally, we say the function $g$ satisfies the quadratic growth condition at $\ox$ if it satisfies this condition at $\ox$ for any $\ov\in \partial g(\ox)$. 
\end{Definition}
As the function $g$ is l.s.c., the set $S(\ox,\ov)$ is closed and $\ox\in S(\ox,\ov)$. When  quadratic growth condition \eqref{LQG} holds, it is clear that 
\begin{equation}\label{Ge}
S(\ox,\ov)\cap\B_\ve(\ox)\subset S\cap \B_\ve(\ox)\quad\mbox{for some}\quad \ve>0.  
\end{equation}
Moreover, for any closed set $S$ fulfilling \eqref{Ge} and $x\in \B^\delta_{\frac{\ve}{2}}(\ox)$, we find some  $u\in S(\ox,\ov)$ such that 
\[
\dist(x;S(\ox,\ov))=\|x-u\|\le \|x-\ox\|< \frac{\ve}{2},
\]
which implies that $\|u-\ox\|\le \|x-\ox\|+\frac{\ve}{2}<\ve$, i.e., $u\in S(\ox,\ov)\cap \B_\ve(\ox)$. It follows from \eqref{Ge} that 
\[
\dist(x;S(\ox,\ov))=\|x-u\|\ge \dist(x;S(\ox,\ov)\cap \B_\ve(\ox))\ge \dist(x;S)
\]
for any $x\in \B^\delta_{\frac{\ve}{2}}(\ox)$. Hence,  if the function $g$ satisfies the quadratic growth condition at $\ox$ for $\ov$, it also satisfies the quadratic growth condition at $\ox$ for $\ov$ w.r.t. any closed set $S$ fulfilling \eqref{Ge}. Many necessary and sufficient conditions for the quadratic growth condition have been established in \cite{BI95a,BS00} and \cite{SW99, W94} under a different name {\em weak sharp minima} with order $2$.

When $g$ is convex and $\ov\in \partial g(\ox)$, the set $S(\ox,\ov)$ coincides with $(\partial g)^{-1}(\ov)=\partial g^*(\ov)$. The quadratic growth condition  of $g$ at $\ox$ to $\ov\in \partial g(\ox)$ w.r.t. $(\partial g)^{-1}(\ov)$  has been studied and connected with the so-called \emph{\L{}ojasiewicz inequality with exponent $\frac{1}{2}$} \cite{BNPS17} and the  \emph{metric subregularity of the subdifferential}   \cite{AG08, DMN14, ZT95} (even for nonconvex cases.) There are  broad classes of convex functions satisfying the quadratic growth condition such as \emph{piece-wise linear quadratic convex} functions \cite[Definition 10.20]{RW98} and many {\em convex spectral} functions \cite{CDZ17}; see also \cite{ZS17} for some other ones.

When $g$ is not convex, the quadratic growth condition  of $g$ at $\ox$ to $\ov\in \partial g(\ox)$ w.r.t. $(\partial g)^{-1}(\ov)$ is the same with the quadratic growth condition of $g$ at $\ox$ for $\ov$ provided that 
\begin{equation}\label{Sep}
(\partial g)^{-1}(\ov)\cap\B_\ve(\ox)\subset S(\ox,\ov)\cap \B_\ve(\ox)\quad \mbox{for sufficiently small}\quad \ve>0.
\end{equation}
It is necessary for the quadratic growth condition \eqref{LQG} at $\ox$ for $\ov$ that $\ox$ is a local minimizer to the function $g_{\ov}(x)\eqdef g(x)-\la \ov,x\ra$, $x\in \XX$.  Then,  condition \eqref{Sep} is similar to the so-called {\em proper separation of isocost surface} of $g_{\ov}$ in \cite{WYYZZ22}, which is an improvement of the {\em proper separation of stationary points} of  $g_{\ov}$ in \cite{LT93}. By \cite[Theorem~3.1]{DMN14}, the function $g$ satisfies the quadratic growth condition  at $\ox$ for $\ov\in \partial g(\ox)$ w.r.t. $(\partial g)^{-1}(\ov)$ provided that $\partial g$ is {\em metrically subregular} at $\ox$ for $\ov$ in the sense that there exist $\eta,\ell>0$ such that
\begin{equation}\label{MSR}
\dist(x;(\partial g)^{-1}(\ov))\le \ell \dist(\ov;\partial g(x))\quad \mbox{for}\quad x\in \B_\eta(\ox).  
\end{equation}
This condition is satisfied when $\partial g$ is a {\em piecewise polyhedral} set-valued mapping, i.e., the graph of $\partial g$, $\{(x,v)\in \XX\times\XX^*|\;v\in \partial g(x)\}$ is  a union of finitely many polyhedral sets; see, e.g., \cite[Example 9.57]{RW98}. Thus the class of (possibly nonconvex)   piecewise linear-quadratic functions fulfills \eqref{MSR}; see also \cite{WYYZZ22} for several sufficient conditions for \eqref{MSR} and some special nonconvex  piecewise linear-quadratic regularizers such as SCAD and MCD penalty functions. Although our theory in this section is applicable to nonconvex functions, we focus our later applications on  low-rank minimization problems \eqref{LSN} when $g$ is the nuclear norm, which also satisfies the quadratic growth condition \cite{ZS17} but the graph of $\partial g$ is not piecewise polyhedral.





The following lemma plays an important role in our analysis. 

\begin{Lemma}[Necessary condition for quadratic growth]\label{MLem} Let  $g:\XX\to \oR$ be a proper l.s.c. function and $S$ be a closed subset of $\XX$ with $\ox\in \dom g\cap S$. If  $g$ satisfies the quadratic growth at $\ox$ for some $\ov\in \partial g(\ox)$ w.r.t. $S$ with some modulus $\kk>0$, we have 
\begin{equation}\label{S1}
    d^2g(\ox|\ov)(w)\ge \kk [\dist(w;T_S(\ox))]^2\quad\mbox{for all}\quad w\in \XX.
    \end{equation}
Moreover, if $g$ satisfies the quadratic growth at $\ox$ for  $\ov$, we have
\begin{equation}\label{Dom}
\Ker d^2 g(\ox|\ov)\eqdef\left\{w\in \XX|\;d^2 g(\ox|\ov)(w)=0\right\}=T_{S(\ox,\ov)}(\ox).
\end{equation}
\end{Lemma}
\noindent{\bf Proof.} Suppose that inequality \eqref{LQG} holds with some $\ve,\delta,\kk>0$. Pick $w\in \XX$, we only need to verify \eqref{S1} when $d^2 g(\ox|\ov)(w)<\infty$, i.e., $w\in \dom d^2 g(\ox|\ov)$. It follows from \eqref{Subd} that there exist sequences $t_k\dn 0$ and $w_k\to w$ such that 
\begin{equation}\label{LimSD}
 d^2 g(\ox|\ov)(w)=\lim_{k\to \infty}\dfrac{g(\ox+t_kw_k)-g(\ox)-t_k\la \ov, w_k\ra}{\frac{1}{2}t_k^2}.
\end{equation}
Hence, we have
\[
g(\ox+t_kw_k)-g(\ox)-\la \ov, \ox+t_kw_k-\ox\ra<\delta
\]
 when $k$ is sufficiently large. Combining \eqref{LQG} and \eqref{LimSD} gives us that 
 \begin{eqnarray}\label{in1}\begin{array}{ll}
    d^2 g(\ox|\ov)(w)&\disp\ge  \kk\limsup_{k\to \infty}\left[\dfrac{\dist(\ox+t_kw_k;S)}{t_k}\right]^2\\
    &\disp=  \kk\limsup_{k\to \infty}\left[\dist\left(w_k;\frac{S-\ox}{t_k}\right)\right]^2
    \end{array}
\end{eqnarray}
As $S$ is closed, there exist $u_k\in \dfrac{S-\ox}{t_k}$, i.e., $\ox+t_ku_k\in S$ such that $
\dist\left(w_k;\dfrac{S-\ox}{t_k}\right)= \|w_k-u_k\|.$ 
This together with \eqref{in1} implies that
\[
d^2 g(\ox|\ov)(w)\ge\kk(\limsup_{k\to \infty} \|w_k-u_k\|^2).
\]
Hence $u_k$ is bounded. By passing to a subsequence, we suppose that $u_k$ converges to $u\in \XX$. As $\ox+t_ku_k\in S$, we have $u\in T_S(\ox)$. It follows from the above inequality that
\[
d^2 g(\ox|\ov)(w)\ge \kk\|w-u\|^2\ge \kk[\dist(w;T_S(\ox))]^2, 
\]
which verifies \eqref{S1}. 

To justify \eqref{Dom}, suppose further that $g$ satisfies the quadratic growth condition at $\ox$ for $\ov$. For any $w\in \Ker  d^2 g(\ox|\ov)$,  we obtain from 
\eqref{S1} that $\dist(w;T_{S(\ox,\ov)}(\ox))=0$, which means $w\in T_{S(\ox,\ov)}(\ox)$. It follows that $\Ker d^2 g(\ox|\ov)\subset T_{S(\ox,\ov)}(\ox)$.   Let us prove the opposite inclusion by picking any $w\in T_{S(\ox,\ov)}(\ox)$. There exist $t_k\dn 0$ and $w_k \to w$ such that $\ox+t_kw_k\in S(\ox,\ov)$, i.e., 
\[
g(\ox+t_kw_k)-g(\ox)-t_k\la \ov,w_k\ra\le 0
\]
We obtain from \eqref{S1} that  
\begin{equation}\label{Zero}
0\le\kk [\dist(w;T_{S(\ox,\ov)}(\ox))]^2\le d^2 g(\ox|\ov)(w)\le \liminf_{k\to \infty} \dfrac{g(\ox+t_kw_k)-g(\ox)-t_k\la \ov,w_k\ra}{\frac{1}{2}t_k^2}\le 0,
\end{equation}
which yields $w\in \Ker d^2 g(\ox|\ov)$ and  verifies $T_{S(\ox,\ov)}(\ox)\subset \Ker d^2 g(\ox|\ov)$. The proof is complete. \endproof

Next let us establish the main theorem of this section, which provides a geometric characterization for strong minima of  problem \eqref{P}. 

\begin{Theorem}[Necessary and sufficient conditions for strong minima]\label{MTheo} Let  $\ox\in {\rm int}\,(\dom f)\cap \dom g$ be a stationary point of  problem \eqref{P}.  If $\ox$ is a strong solution of  problem \eqref{P},  then 
\begin{equation}\label{SC}
\la \nabla^2 f(\ox)w,w\ra>0\qquad \mbox{for all}\quad w\in T_{S(\ox,\ov)}(\ox)\setminus\{0\}. 
\end{equation}
Suppose further that $g$ satisfies the quadratic growth  condition at $\ox$ for $\ov\eqdef-\nabla f(\ox)$ and that $\nabla^2 f(\ox)$ is positive semidefinite, then $\ox$ is a strong solution of problem \eqref{P} if and only if  
\begin{equation}\label{SC2}
\Ker \nabla^2 f(\ox)\cap T_{S(\ox,\ov)}(\ox)=\{0\}. 
\end{equation}
\end{Theorem}
\noindent{\bf Proof.} As $f$ is twice continuously differentiable at $\ox\in {\rm int}\, (\dom f)$, we derive from    \eqref{SS} and \eqref{SR2} that $\ox$ is a strong solution of $\ph$ if and only if there exists some $\ell>0$ such that
\begin{equation}\label{SQ}
\la \nabla^2 f(\ox)w,w\ra+ d^2g(\ox|\ov)(w)\ge  \ell\|w\|^2 \quad \mbox{for all}\quad w\in \XX. 
\end{equation}
 To justify the first part, suppose that $\ox$ is a strong solution of $\ph$, i.e., \eqref{SQ} holds. Pick any $w\in  T_{S(\ox,\ov)}(\ox)\setminus\{0\}$ and find sequences $t_k\dn 0$ and $w_k\to w$ such that $\ox+t_kw_k\in  S(\ox,\ov)$, which means
\[
\dfrac{g(\ox+t_kw_k)-g(\ox)-\la \ov, \ox+t_kw_k-\ox\ra}{\frac{1}{2}t_k^2}\le 0
\]
By the definition of $d^2 g(\ox|\ov)(w)$ in \eqref{Subd}, we have $d^2 g(\ox|\ov)(w)\le 0$.
 This together with \eqref{SQ} verifies \eqref{SC}.

To verify the second part of the theorem, suppose that the function $g$ satisfies the quadratic growth condition at $\ox$ for $\ov$ with some modulus $\kk>0$ and $\nabla f^2(\ox)$ is positive semidefinite. It is obvious that \eqref{SC} implies \eqref{SC2}. We only need to prove that \eqref{SC2} is sufficient for strong minima at $\ox$. Suppose that condition \eqref{SC2} is satisfied. 
If  condition \eqref{SQ} failed, we could find  a sequence $w_k$ such that $\|w_k\|=1$ and   
\[
\la \nabla^2 f(\ox)w_k,w_k\ra+ d^2g(\ox|\ov)(w_k)\leq \frac{1}{k}.
\]
It follows from \eqref{S1} that 
\begin{equation}\label{in2}
\frac{1}{k}\ge \la \nabla^2 f(\ox)w_k,w_k\ra+\kk [\dist(w_0;T_{S(\ox,\ov)}(\ox))]^2.
\end{equation}
By passing to a subsequence, assume that $w_k\to w_0$ with $\|w_0\|=1$ (without relabeling.) It follows that
\[
0\ge \la \nabla^2 f(\ox)w_0,w_0\ra+\kk [\dist(w_0;T_{S(\ox,\ov)}(\ox)]^2\ge \la \nabla^2 f(\ox)w_0,w_0\ra\ge 0.
\]
Hence, we have $\la \nabla^2 f(\ox)w_0,w_0\ra=0$ and $\dist(w_0;T_{S(\ox,\ov)}(\ox))=0$, which means 
\[
w_0\in \Ker \nabla^2 f(\ox) \cap  T_{S(\ox,\ov)}(\ox).
\]
This is a contradiction to \eqref{SC2} as $\|w_0\|=1$. Hence, \eqref{SC2} holds for some $\ell>0$ and $\ox$ is a strong solution of $\ph$. The proof is complete. \endproof

\begin{Corollary}[Geometric characterization for strong minima of  convex problems]\label{MCoro} Let $f,g:\XX \to \oR$ be proper l.s.c. convex functions and $\ox\in {\rm int}\, (\dom f)\cap \dom g$ be a minimizer to problem \eqref{P}. Suppose that $f$ is twice continuously differentiable in ${\rm int}\, (\dom f)$  and that $\partial g$ is metrically subregular at $\ox$ for $\ov=-\nabla f(\ox)$.  Then  $\ox$ is a strong solution to problem \eqref{P} if any only if 
\begin{equation}\label{SC3}
\Ker \nabla^2 f(\ox)\cap T_{(\partial g^*)(\ov)}(\ox)=\{0\}. 
\end{equation}
\end{Corollary}
{\noindent \bf Proof.} As discussed before \eqref{MSR}, when $g$ is a convex function,  the metric subregularity of $\partial g$ at $\ox$ for $\ov$ implies the quadratic growth condition (they are indeed equivalent \cite{AG08,ZT95}.) 
Since $f$ is convex, $\nabla^2 f(\ox)$ is positive semidefinite. By Theorem~\ref{MTheo}, $\ox$ is a strong solution if and only if \eqref{SC2} holds. 
Since $g$ is a convex function, we have $S(\ox,\ov)=(\partial g)^{-1}(\ov)=\partial g^*(\ov)$. Thus \eqref{SC3} is equivalent to \eqref{SC2}. The proof is complete. \endproof

Unlike many other necessary and sufficient conditions for strong minima, our geometric characterizations \eqref{SC2} and \eqref{SC3} do not involve the "curvature" or the "sigma-term" of the function $g$. We still need to compute the contingent cones $T_{S(\ox,\ov)}(\ox)$ in \eqref{SC2} or $T_{\partial g^*(\ov)}(\ox)$ in \eqref{SC3}.   In Section 4 and 5, we consider the case $g=\|\cdot\|_*$, the nuclear norm and  provide a simple calculation of $T_{\partial g^*(\ov)}(\ox)$.    

\vspace{0.1in}

\subsection{Geometric characterization for strong minima of optimization problems with linear constraints}

In this subsection, we apply the idea in Theorem~\ref{MTheo} and Corollary~\ref{MCoro} to the following   convex optimization problem with linear constraints
\begin{equation}\label{CO}
    \min_{x\in \XX} \quad g(x)\quad \mbox{subject to}\quad \Phi x\in K,
\end{equation}
where $g:\mathbb{X}\to \R$ is a continuous (nonsmooth) convex function with full domain, $\Phi:\mathbb{X}\to \mathbb{Y}$ is a linear operator between two Euclidean spaces, and $K$ is a closed convex polyhedral  set in $\mathbb{Y}$. Unlike problem \eqref{P},  function $g$ needs to satisfy more properties such as convexity, quadratic growth condition, and second order regularity stated in the next theorem. 

Let us recall  that  $x_0$ is   a stationary solution of  problem \eqref{CO} if there exists a Lagrange multiplier $\lm\in \mathbb{Y}^*$, hence a \emph{dual certificate}, such that  
\begin{equation}\label{Lag}
-\Phi^*\lm \in \partial g(x_0)\quad\mbox{and}\quad \lm\in N_K(\Phi x_0).
\end{equation}
The set of Lagrange multipliers  is defined by
\begin{eqnarray}\label{Las}
    \Lm(x_0)\eqdef\{\lm\in N_K(\Phi x_0)| -\Phi^*\lm \in \partial g(x_0)\}. 
\end{eqnarray}
The critical cone of  this problem at the  stationary point $x_0$ is 
\begin{equation}\label{Cr}
    C(x_0)\eqdef\{w\in \mathbb{X}|\; \Phi w\in T_{K}(\Phi x_0),  dg(x_0)(w)= 0\}.
\end{equation}
The point $x_0$ is call a strong solutions of problem \eqref{CO} if there exists $\ve>0$ and $c>0$ such that 
\[
g(x)-g(x_0)\ge c\|x-x_0\|^2\qquad \mbox{when}\qquad  \Phi x\in K\quad \mbox{and}\quad x\in \B_\ve(x_0).
\]
\begin{Theorem}[Geometric  characterization for strong minima of problem \eqref{CO}]\label{thm5} Let $x_0$ be a stationary point of  problem \eqref{CO}. Suppose that the convex function $g$ is second order regular at $x_0$ and  satisfies the  quadratic grown condition at $x_0$. Then $x_0$ is a strong solution of  problem \eqref{CO} if and only if  
\begin{equation}\label{CharCo}
      \left[\bigcap_{\lm\in \Lm(x_0)} T_{\partial g^*(-\Phi^*\lm)}(x_0)\right]\cap C(x_0)=\{0\}. 
\end{equation}
\end{Theorem}
\noindent{\bf Proof.} Define $y_0\eqdef\Phi x_0, $ $\mathbb{L}\eqdef\Im\Phi$ being a linear subspace of $\YY$, $K_{\mathbb{L}}\eqdef K\cap \mathbb{L}$,  the mapping $F:\XX\to \mathbb{L}\times \XX$ by $F(x)\eqdef(\Phi x,x)$ for $x\in \XX$, and the function $G:\mathbb{L}\times \XX\to \R$ by $$G(u,x)=\iota_{K_\mathbb{L}}(u)+g(x)\quad \mbox{for}\quad (u,x)\in \mathbb{L}\times \XX.$$ Hence we can replace $\YY$ by $\mathbb{L}$ and $K$ by $K_{\mathbb{L}}$ in problem \eqref{CO}.  Rewrite problem \eqref{CO} as a composite optimization problem \eqref{CoP}
\begin{equation}\label{COP}
    \inf_{x\in \XX}\quad G(F(x)).  
\end{equation}
Observe that  $x_0$ is a strong solution of \eqref{CO} if and only if it is a strong solution  of \eqref{COP}.  Robinson's constraint qualification \eqref{Rob} for problem \eqref{COP} is
\begin{equation}\label{RCQ}
0\in {\rm int}\, (F(x_0)+\nabla F(x_0)\XX-K_\mathbb{L}\times \XX).
\end{equation}
As $\Phi:\XX\to\mathbb{L}$ is a surjective  operator, \cite[Corollary~2.101]{BS00} tells us that  the above condition is equivalent to the existence of $w\in \XX$ satisfying
\begin{equation*}
y_0+\Phi w\in K_\mathbb{L}\quad \mbox{and}\quad x_0+w\in {\rm int}\, \XX=\XX. 
\end{equation*}
This condition holds trivially at $w=0_\XX$. Thus Robinson's constraint qualification  \eqref{RCQ} holds.

The critical cone \eqref{CriC} for problem \eqref{COP} is 
\begin{eqnarray}\label{Crit}
 \Hat C(x_0)\eqdef\{w\in \XX|\;dG(F(x_0)|\; \nabla F(x_0)w)=0\}=\{w\in\XX|\; \Phi w\in T_{K_\mathbb{L}}(y_0), dg(x_0)(w)= 0\}.
\end{eqnarray}
As $K_\mathbb{L}$ is also a  polyhedral in $\mathbb{L}$, we have
\begin{equation}\label{KL1}
T_{K_\mathbb{L}}(y_0)=\R_+(K_{\mathbb L}-y_0)=\R_+(K-y_0)\cap \mathbb{L}=T_K(y_0)\cap \mathbb{L}. 
\end{equation}
It follows that the set $\Hat C(x_0)$ is exactly the critical cone $C(x_0)$ defined in \eqref{Cr}.

By \eqref{Lagr}, the set of Lagrange multipliers of problem \eqref{COP} is 
\begin{eqnarray}\begin{array}{ll}\label{Lm}
\Hat \Lm(x_0)&\eqdef\{(\lm,\mu)\in \mathbb{L}^*\times \XX^*|\; \nabla F(x_0)^*(\lm,\mu)=0,(\lm,\mu)\in \partial G(F(x_0))\}\\
&=\{(\lm,\mu)\in \mathbb{L}^*\times \XX^*|\; \mu=-\Phi^*\lm,\lm\in N_{K_{\mathbb{L}}}(y_0), \mu\in \partial g(x_0)\}.
\end{array}
\end{eqnarray}
Note further that 
\[
\epi G=K_\mathbb{L}\times \epi g.
\]
Since $\mathbb{L}$ is a polyhedral, it is second order regular \cite{BCS99}. As $\epi g$ is second order regular, so is $\epi G$; see, e.g., \cite[Propisition 3.89]{BS00}.

By Theorem~\ref{ComP}, $x_0$ is a strong solution of problem \eqref{COP} if and only if for any $w\in C(x_0)\setminus\{0\}$ there exists $(\lm,\mu)\in \Lm(x_0)$ such that 
\begin{eqnarray}\label{FPsi}
\la (\lm,\mu),\nabla^2 F(x_0)(w,w)\ra+d^2G(F(x_0)|(\lm,\mu))(\nabla F(x_0)w)>0.
\end{eqnarray}
Observe that 
\begin{equation}\label{G2}
d^2G(F(x_0)|\nabla F(x_0)(\lm,\mu))(\nabla F(x_0)w)=d^2g(x_0|\mu)(w)+d^2\iota_{K_\mathbb{L}}(y_0|\lm)(\Phi w).
\end{equation}
Note from \eqref{Subd} that   
\begin{equation}\label{TKL}
    d^2\iota_{K_\mathbb{L}}(y_0|\lm)(\Phi w)=\liminf_{z \to \Phi w, t\dn 0}\frac{\iota_{K_\mathbb{L}}(y_0+tz)-\iota_{K_\mathbb{L}}(y_0)- t\la \lm,z\ra}{0.5t^2}\ge 0. 
\end{equation}
By \eqref{KL1}, we have
\begin{equation}\label{TKL2}
    d^2\iota_{K_\mathbb{L}}(y_0|\lm)(\Phi w)=\liminf_{z \st{T_{K_\mathbb{L}}(y_0)}\to \Phi w, t\dn 0}\frac{-\la \lm,z\ra}{0.5t}\ge 0. 
\end{equation}
 Since $w\in C(x_0)$, we have $\Phi w\in T_{K_\mathbb{L}}(\Phi x_0)$. As $\lm \in N_{K_\mathbb{L}}(\Phi x_0)$ and $\mu=-\Phi^*\lm\in \partial g(x_0)$, it follows that 
\[
0= dg(x_0)(w)\ge \la \mu, w\ra=-\la \Phi^*\lm ,w\ra=-\la \lm,\Phi w\ra\ge 0,
\]
which implies that $\la \lm, \Phi w\ra=0$. This together with \eqref{TKL2} tells us that $d^2\iota_{K_\mathbb{L}}(y_0|\lm)(\Phi w)=0$. 

By \eqref{G2}, condition \eqref{FPsi} is equivalent to 
\begin{equation}\label{gS}
d^2g(x_0|-\Phi^*\lm)(w)>0. 
\end{equation}
Since $K$ is a polyhedral set, we have 
\[
N_{K_\mathbb{L}}(\Phi x_0)=N_{K\cap \mathbb{L}}(\Phi x_0)=N_K(\Phi x_0)+N_\mathbb{L}(\Phi x_0)=N_K(\Phi x_0)+\Ker  \Phi^*.
\]
Represent $\lm=\lm_1+\lm_2$ for some $\lm_1\in N_K(\Phi x_0)$ and $\lm_2\in \Ker  \Phi^*$, it follows that
\[
-\Phi^*\lm=-\Phi^*(\lm_1+\lm_2)=-\Phi^*\lm_1.
\]
Hence $x_0$ is a strong solution of problem \eqref{COP} if any only if for any $w\in C(x_0)\setminus \{0\}$ there exists some $
\lm\in \Lm(x_0)$ in \eqref{Lag} such that \eqref{gS} holds. Since $g$ satisfies the quadratic grown condition at $x_0$, we get from Lemma~\ref{MLem} that 
\[
\Ker  d^2g(x_0|-\Phi^*\lm)=T_{\partial g^*(-\Phi^*\lm)}(x_0).
\]
Hence condition \eqref{gS} means  for any $w\in C(x_0)\setminus\{0\}$ there exists $\lm\in \Lm(x_0)$ such that 
\[
w\notin T_{\partial g^*(-\Phi^*\lm)}(x_0)\quad \mbox{or}\quad w\in \XX\setminus T_{\partial g^*(-\Phi^*\lm)}(x_0). 
\]
This is equivalent to the following inclusion
\begin{equation*}
   C(x_0)\setminus\{0\}\subset  \left[\bigcup_{\lm\in \Lm(x_0)} \left(\XX\setminus T_{\partial g^*(-\Phi^*\lm)}(x_0)\right) \right]=\XX\setminus \left[\bigcap_{\lm\in \Lm(x_0)}T_{\partial g^*(-\Phi^*\lm)}(x_0)\right], 
\end{equation*}
which is also equivalent to \eqref{CharCo}. The proof is complete. \endproof 

The approach of using composite function \eqref{COP} to study problem \eqref{CO} is traditional; see, e.g., \cite{BCS99, I79}. However, by assuming additionally the function $g$ satisfies the quadratic grown condition at  $x_0$, we are able to obtain the new geometric characterization for strong solution in \eqref{CharCo}. The main  idea in this result  is similar to that in Theorem~\ref{MTheo}. Here we require the function $g$ to satisfy more assumptions, but when applying this result to the nuclear norm minimization problem \eqref{NNM}, they are also valid, as the nuclear norm is second order regular and also satisfies the quadratic growth condition \cite{CDZ17}.

\section{Characterizations for strong minima of  low-rank optimization problems}
\setcounter{equation}{0}
This section devotes to new  characterizations for strong minima of the  low-rank optimization problem:
\begin{equation}\label{Nucl}
    \min_{X\in \R^{n_1\times n_2}}\quad h(\Phi X)+\mu \|X\|_*,
\end{equation}
where $\Phi:\R^{n_1\times n_2}\to \R^{m}$ is a linear operator,  $g(X)\eqdef\|X\|_*$ is the nuclear norm of $X\in \R^{n_1\times n_2}$, $\mu$ is a positive constant,  and $h:\R^m\to \oR$ satisfies the following standing assumptions \cite{FNT21}:
\begin{enumerate}[{\rm (A)}]
    \item $h$ is proper convex and twice continuously  differentiable in ${\rm int}\,(\dom h)$.
    \item $\nabla^2 h(\Phi X)$ is positive definite for any $X\in \Phi^{-1}({\rm int}\,(\dom h))$. 
\end{enumerate}
Strongly convex functions with full domain clearly satisfy the above standing assumptions. Another important (non-strongly convex)  function with these conditions widely used in statistical/machine learning is the {\em Kullback-Leiber divergence}. Sufficient conditions for strong minima of   problem \eqref{Nucl} can be obtained from \cite[Theorem~12]{CDZ17}. However, their result still relies on some  computation  of $d^2\|\cdot\|_*$, which is  complicated; see, e.g.,  \cite{D17,ZZX13} and the recent paper \cite{MS23} for the case of symmetric matrices. We will provide some explicit and computable characterizations for strong minima of problem \eqref{Nucl}  based on  Corollary~\ref{MCoro}.   The calculation of the contingent cone $T_{\partial g^*(\OY)}(\OX)$ is rather simple; see our formula \eqref{Tang2} below. 

Let us recall a few standard notations for matrices. The space of all matrices $\R^{n_1\times n_2}$ ($n_1\le n_2$) is endowed with the inner product
\begin{equation*}
    \la X,Y\ra\eqdef\Tr(X^TY)\quad \mbox{for all}\quad X,Y\in \R^{n_1\times n_2},
\end{equation*}
where $\Tr$ is the \emph{trace operator}. The  Frobenious norm on $\R^{n_1\times n_2}$ is 
\[
\|X\|_F\eqdef\sqrt{\Tr(X^TX)}\quad \mbox{for all}\quad X\in \R^{n_1\times n_2}. 
\]
The nuclear norm and spectral norm of $X\in \R^{n_1\times n_2}$ are defined respectively by
\[
\|X\|_*\eqdef\sum_{i=1}^{n_1}\sigma_i(X)\quad\mbox{and}\quad  \|X\|\eqdef\sigma_1(X),
\]
where $\sigma_1(X)\ge \sigma_2(X)\ge \ldots\ge \sigma_{n_1}(X)\ge 0$ are all singular values of $X$. 
   Suppose that a {\em full} Singular Value Decomposition (SVD) of $\Bar X\in \R^{n_1\times n_2}$ is 
\begin{equation}\label{OX}
\Bar X=U\begin{pmatrix}\Bar \Sigma_r &0\\
0& 0\end{pmatrix}_{n_1\times n_2}V^T \quad \mbox{with}\quad
\Bar\Sigma_r=\begin{pmatrix}\sigma_1(\OX)& \dots &0\\
\vdots&\ddots &\vdots\\
0 &\ldots & \sigma_r(\OX)\end{pmatrix},
\end{equation}
where $r={\rm rank}\, (\Bar X)$
, $U\in \R^{n_1\times n_1}$ and $V\in \R^{n_2\times n_2}$ are orthogonal matrices. Let $\mathcal{O}(\OX)$ be the set of all such pairs $(U,V)$ satisfying \eqref{OX}. We write $U=\begin{pmatrix} U_I&U_J\end{pmatrix}$ and $V=\begin{pmatrix} V_I&V_K\end{pmatrix}$, where  $U_I$ and $V_I$ are the submatrices of the first $r$ columns of  $U$ and $V$, respectively. We get from \eqref{OX} that $\OX=U_I\Bar\Sigma_rV_I^T$, which is known as a {\em compact SVD} of $\OX$.  

The following lemma is significant in our paper. The first part is well-known \cite[Example~2]{W92}. The last part was established in \cite[Proposition 10]{ZS17}, which can be viewed as a direct consequence of  \cite[Example~1]{W92} via convex analysis, the formula of normal cone to a level set \cite[Corollary 23.7.1]{R70}.

\begin{Lemma}[Subdifferential of  the nuclear norm]\label{Lem} The subdifferential to nuclear norm at $\OX\in \R^{n_1\times n_2}$ is computed by
\begin{equation}\label{subdif}
    \partial\|\OX\|_*=\left\{U\begin{pmatrix} \Id_r&0\\ 0 & W\end{pmatrix}V^T|\; \|W\|\le 1\right\} \quad \mbox{for any}\quad (U,V)\in \mathcal{O}(\OX). 
\end{equation}
Moreover, $\OY\in \partial\|\OX\|_*$ if and only if $\|\OY\|\le 1$ and 
\begin{equation}\label{Fen}
    \|\OX\|_*=\la \OY,\OX\ra. 
\end{equation}
Furthermore, for any $\OY\in\B\eqdef\{Z\in \R^{n_1\times n_2}|\;\|Z\|\le 1\}$, we have 
\begin{equation}\label{Inver}
\partial g^*(\OY)= N_\B(\OY)=\OU\begin{pmatrix}\mathbb{S}_+^{p(\OY)} &0\\0 &0\end{pmatrix}\OV^T\quad \mbox{for any}\quad (\OU,\OV)\in \mathcal{O}(\OY),
\end{equation}
where $\mathbb{S}_+^p$ is the set of all $p\times p$  symmetric positive semidefinite matrices and  $p(\OY)$ is defined by
\begin{equation}\label{p}
p(\OY)\eqdef\#\{i|\; \sigma_i(\OY)=1\}.
\end{equation}
\end{Lemma}

  Let $\Bar Y\in \partial \|\Bar X\|_*$ and $(U,V)\in \mathcal{O}(\OX)$. It follows from \eqref{subdif} that $\OY$  can be represented by
\begin{equation}\label{OY}
\Bar Y=U\begin{pmatrix}\Id_r &0\\
0& \Bar W\end{pmatrix}V^T
\end{equation}
with some $\Bar W\in \R^{(n_1-r)\times (n_2-r)}$ satisfying $\|\Bar W\|\le 1$. 
Let $(\Hat U,\Hat V)\in \mathcal{O}(\Bar W)$ and $\Hat U\Sigma \Hat V^T$ be a full SVD of $\OW$.  We get from \eqref{OY} that 
\begin{equation}\label{OY2}
\Bar Y=\OU\begin{pmatrix}\Id_r &0\\
0& \Sigma\end{pmatrix}\OV^T\quad \mbox{with}\quad \OU\eqdef(U_I\; U_J\Hat U)\quad \mbox{and}\quad \OV\eqdef(V_I \; V_K\Hat V).
\end{equation}
Observe that  $\OU^T\OU=\Id_{{n_1}}$ and $\OV^T\OV=\Id_{n_2}$. It follows that $(\OU,\OV)\in \mathcal{O}(\OX)\cap \mathcal{O}(\OY)$, which means $\OX$ and $\OY$ have {\em simultaneous ordered singular value decomposition} \cite{LS05a,LS05b} with orthogonal matrix pair $(\OU,\OV)$ in the sense that 
\begin{equation}\label{SSVD}
    \OX=\OU({\rm Diag}\,\sigma(\OX))\OV^T\qquad \mbox{and}\qquad \OY=\OU({\rm Diag}\,\sigma(\OY))\OV^T,
\end{equation}
where $\sigma(\OX)\eqdef\left (\sigma_1(\OX), \ldots,\sigma_{n_1}(\OX)\right)^T$ and ${\rm Diag}\,\sigma(\OX)\eqdef\begin{pmatrix}\sigma_1(\OX)&\ldots&0&0&\ldots& 0\\
0&\ddots&0&0&\ldots& 0\\ 
0&\ldots &\sigma_{n_1}(\OX)&0&\ldots 
&0\end{pmatrix}_{n_1\times n_2}$.

The following result establishes a geometric characterization for strong solution of the problem \eqref{Nucl}. According to \cite[Proposition~11]{ZS17}, the subdifferential of the  nuclear norm function  satisfies the {\em metric subregularity} \eqref{MSR} at any $\OX$ for any $Y\in \partial\|\OX\|_*$.

\begin{Corollary}[Geometric characterization for strong minima of  low-rank optimization problems]\label{CoNu} Suppose that  $\Bar X\in \Phi^{-1}({\rm int}(\dom h))$ is a minimizer of  problem \eqref{Nucl} with $\Bar Y\eqdef-\frac{1}{\mu}\Phi^*\nabla h(\Phi\Bar X)\in \partial \|\Bar X\|_*$. Let   $(\OU,\OV)\in \mathcal{O}(\OX)\cap\mathcal{O}(\OY)$ as in \eqref{SSVD} or \eqref{OY2}.   Then we have 
\begin{equation}\label{Tang2}
T_{N_\B(\Bar Y)}(\Bar X)=\left\{\OU\begin{pmatrix}A& B &0\\ B^T&C& 0\\0&0&0\end{pmatrix}\OV^T|\; A\in \mathbb{S}^{r}, B\in \R^{r\times (p(\OY)-r)},  C\in \mathbb{S}^{p(\OY)-r}_+\right\},
\end{equation}
where $p(\OY)$ is defined in \eqref{p} and $\mathbb{S}^r$ is the set of all symmetric matrices of size $r\times r$.
Hence $\Bar X$ is a strong solution of  \eqref{Nucl} if any only if 
\begin{equation}\label{Char1}
\Ker \Phi \cap T_{N_\B(\Bar Y)}(\Bar X)=\{0\}. 
\end{equation}
Consequently, $\OX$ is a strong solution of  \eqref{Nucl} provided that the following \emph{Strong Sufficient Condition} holds
\begin{equation}\label{Suff}
\Ker\Phi\cap \OU\begin{pmatrix}\mathbb{S}^{p(\OY)}&0\\0&0\end{pmatrix}\OV^T=\{0\}. 
\end{equation}
\end{Corollary}
\noindent{\bf Proof.} By Lemma~\ref{Lem}, we have 
\[
\partial g^*(\OY)=N_\B(\OY)=\OU\begin{pmatrix}\mathbb{S}_+^{p(\OY)} &0\\0 &0\end{pmatrix}\OV^T.
\]
As $(\OU,\OV)\in \mathcal{O}(\OX)$, we obtain from \eqref{OX} that
\[
T_{N_\B(\OY)}(\OX)=\OU \begin{pmatrix}T_{\mathbb{S}_+^{p(\OY)}}\begin{pmatrix}\Bar\Sigma_r &0\\0 &0\end{pmatrix} &0\\0 &0\end{pmatrix}\OV^T,
\]
which is exactly the right-hand side of \eqref{Tang2} according to the contingent cone formula to $S^{p(\OY)}_+$ in \cite[Example~2.65]{BS00}.  Since $\partial \|\cdot\|_*$ is metrically subregular at $\OX$ for $\OY$ by \cite[Proposition~11]{ZS17}, it follows from  Corollary~\ref{MCoro} that $\OX$ is a strong solution of  problem \eqref{Nucl} if and only if 
\[
\Ker\Big(\Phi^*\nabla^2 h(\Phi \OX)\Phi\Big) \cap T_{N_\B(\OY)}(\OX)=\{0\}.
\]
Since $\nabla^2 h(\Phi \OX)\succ0$, we have $\Ker\Big(\Phi^*\nabla^2 h(\Phi \OX)\Phi\Big)=\Ker \Phi$. The characterization \eqref{Char1} for strong minima at $\OX$ follows from the above condition.

Finally, note from \eqref{Tang2} that 
\[
T_{N_\B(\OY)}(\OX)\subset \OU\begin{pmatrix}\mathbb{S}^{p(\OY)}&0\\0&0\end{pmatrix}\OV^T.
\]
Strong Sufficient Condition \eqref{Suff} implies strong minima at $\OX$ by \eqref{Char1}.  \endproof

\begin{Remark}{\rm The geometric characterization \eqref{Char1} for strong minima of problem \eqref{Nucl} is news. A sufficient condition is indeed obtained from \cite[Theorem~12]{CDZ17}, which considers more general optimization problems involving {\em spectral functions}. However, their result contains a nontrivial {\em sigma-term}, which is calculated explicitly in recent papers \cite{CD19,MS23} for the case of symmetric matrices.  Our approach is totally different without any  sigma-terms. Moreover,  our condition is a full  characterization for strong minima. 

Another result about strong minima of problem \eqref{Nucl} was established in  \cite[Proposition~12]{LFP17}, which plays an important role in proving the local linear convergence of Forward-Backward algorithms solving problem \eqref{Nucl}. The result mainly states that  the so-called {\em Restricted Injectivity} and {\em Nondegenerate Condition} are sufficient for strong minima; see also \cite[Proposition~4.27]{FNT21} for similar observation.  Let us recall these important conditions here; see further discussions about them in Section~5. Let $(U,V) \in \mathcal{O}(\OX)$ and  define the {\em model tangent subspace} 
\begin{equation}\label{T}
    \TT\eqdef\{U_I Y^T+XV_I^T|\; X\in \R^{n_1\times r}, Y\in \R^{n_2\times r}\}
\end{equation} of $\R^{n_1\times n_2}$ with dimension $\dim \TT=r(n_1+n_2-r)$; see, e.g.,   \cite{CR09,CR13}. The Restricted Injectivity condition means 
\begin{equation}\label{RI1}
\Ker \Phi\cap \TT=\{0\}. 
\end{equation}
And the Nondegeneracy Condition holds when  
\begin{equation}\label{NC1}
    \OY=-\frac{1}{\mu}\Phi^*\nabla h(\Phi \OX)\in {\rm ri}\, \partial \|\OX\|_*, 
\end{equation}
where ${\rm ri}\, \partial \|\OX\|_*$ is the \emph{relative interior} of $\partial \|\OX\|_*$; see \cite{R70}. 
The validity of  Nondegeneracy Condition \eqref{NC1} implies that $\OX$ is an optimal solution of problem \eqref{Nucl}. Note that 
\begin{equation}\label{Rei}
{\rm ri}\, \partial \|\OX\|_*=\left\{U\begin{pmatrix}\Id_r&0\\0&W\end{pmatrix}V^T|\; \|W\|<1\right\}\quad \mbox{with}\quad r=\rank(\OX). 
\end{equation}
Hence, Nondegeneracy Condition \eqref{NC1} means that  the number of singular value ones, $p(\OY)$ in \eqref{p}, is the rank of $\OX$. In this case,  Restricted Injectivity \eqref{RI1} clearly implies  the Strong Sufficient Condition \eqref{Suff}. Hence, the combination of Restricted Injectivity \eqref{RI1} and Nondegeneracy Condition \eqref{NC1} is stronger than our Strong Sufficient Condition \eqref{Suff}. The following result gives a complete picture about strong minima when Nondegeneracy Condition \eqref{NC1} occurs.
}
\end{Remark}

\begin{Corollary}[Strong minima under Nondegeneracy Condition]\label{Strict1} Suppose that $\Bar X\in \Phi^{-1}({\rm int}(\dom h))$ and Nondegeneracy Condition \eqref{NC1} holds. Then $\OX$ is a strong solution of problem \eqref{Nucl} if and only if the following {\em Strict Restricted Injectivity} holds
\begin{equation}\label{SRI1}
    \Ker \Phi\cap U_I\mathbb{S}^rV_I^T=\{0\},
\end{equation}
where $U_I\Bar\Sigma V_I^T$ is a compact SVD of $\OX$.
\end{Corollary}
\noindent{\bf Proof.} As Nondegeneracy Condition \eqref{NC1} holds, $\OX$ is a solution of problem \eqref{Nucl}. In this case, observe from \eqref{OY2} and \eqref{Tang2} that $T_{N_\B(\OY)}(\OX)=U_I\mathbb{S}^r V_I^T.$ The equivalence between strong minima at $\OX$ and \eqref{SRI1} follows Corollary~\ref{MCoro}.\endproof

As the dimension of subspace $U_I\mathbb{S}^rV_I^T$ is $\frac{1}{2}r(r+1)$, which is usually small in low-rank optimization problems, it is likely that condition \eqref{SRI1} holds when Nondegeneracy Condition \eqref{NC1} is satisfied. More discussions about Strict Restricted Injectivity will be added on Section~5. 

Although geometric characterization \eqref{Char1} looks simple, checking it in high dimension is  nontrivial. But Strong Sufficient Condition \eqref{Suff} and Strict Restricted Injectivity \eqref{SRI1}  can be verified easily.  Next we  establish some quantitative characterizations for strong minima. Before doing so, we obtain some projection formulas onto subspaces $\TT$ and $\TT^\perp$. For any $X\in \R^{n_1\times n_2}$, suppose that $X$ is represented by block matrices as: 
\begin{eqnarray*}\begin{array}{ll}
X&=\begin{pmatrix}U_I & U_J\end{pmatrix}\begin{pmatrix} A& B\\C&D\end{pmatrix}\begin{pmatrix}V_I & V_K\end{pmatrix}^T\quad \mbox{with}\quad (U,V)\in \mathcal{O}(\OX).
\end{array}
\end{eqnarray*}
The projections of $X$ onto $\TT$ and $\TT^\perp$ are computed respectively by 
\begin{equation}\label{Proj}
P_{\TT}X=   \begin{pmatrix}U_I & U_J\end{pmatrix}\begin{pmatrix} A& B\\C&0\end{pmatrix}\begin{pmatrix}V_I & V_K\end{pmatrix}^T\quad \mbox{and}\quad P_{\TT^\perp}X=U_JDV_K^T.
\end{equation}

The following result provides a formula for critical cone of nuclear norm at $\OX$ for $\OY\in \partial\|\OX\|_*$.

\begin{Proposition}[Critical cone of nuclear norm]\label{Pro4} Let $\Bar Y\in \partial \|\Bar X\|_*$ and $(\OU,\OV)\in \mathcal{O}(\Bar X)\cap \mathcal{O}(\OY)$ as in  \eqref{OX} and \eqref{OY2}. Define $H\eqdef\{k\in \{r+1, \ldots n_1\}|\; \sigma_k(\OY)=1\}$.  Then the critical cone $\mathcal{C}(\OX,\OY)$ of $\|\cdot\|_*$ at $\OX$ for $\OY$ is computed by
\begin{equation}\label{CritY}
\mathcal{C}(\OX,\OY)\eqdef\{W\in \R^{n_1\times n_2}|\; d\|\OX\|_*(W)=\la \Bar Y, W\ra\}=\left\{W\in \R^{n_1\times n_2}|\; P_{\TT^\perp}W\in \OU_H\mathbb{S}^{|H|}_+\OV_H^T\right\},
\end{equation}
where $\OU_H$ and $\OV_H$ are submatrices of index columns $H$ of $\OU$ and $\OV$, respectively.
\end{Proposition}
\noindent{\bf Proof.} For any $W\in \R^{m\times n}$, it is well-known from convex analysis \cite{R70} that  
\begin{eqnarray}\label{ConV}
   d\|\OX\|_*(W)=\sup_{Y\in \partial  \|\OX\|_*}\la Y,W\ra. 
\end{eqnarray}
This together with  \eqref{subdif} and \eqref{Proj} implies that 
\begin{eqnarray}\label{CXY}\begin{array}{ll}
d\|\OX\|_*(W)&\disp=\sup_{Y\in \partial  \|\OX\|_*}\la Y,P_{\TT}W+P_{\TT^\perp}W\ra=\sup_{Y\in \partial  \|\OX\|_*} \la P_{\TT}Y,W\ra+\la Y,P_{\TT^\perp}W\ra\\
&=\la E,W\ra+\|P_{\TT^\perp}W\|_*
\end{array}
\end{eqnarray}
with $E\eqdef U_IV_I^T$. As $\OY\in \partial \|\OX\|_*$,  we have $W\in \mathcal{C}(\OX,\OY)$ if and only if 
\[
\la E,W\ra+\|P_{\TT^\perp}W\|_*=\la P_{\TT}\OY,W\ra+\la P_{\TT^\perp} \OY,W\ra=\la E,W\ra+\la P_{\TT^\perp} \OY,P_{\TT^\perp}W\ra,
\]
which means  $\|P_{\TT^\perp}W\|_*=\la P_{\TT^\perp} \OY,P_{\TT^\perp}W\ra$. By Lemma~\ref{Lem}, we have  $P_{\TT^\perp}W\in \partial  g^*(P_{\TT^\perp} \OY)$, or equivalently  $P_{\TT^\perp}W\in \OU_H\mathbb{S}^{|H|}_+\OV_H^T$. The proof is complete. 
\endproof

Next, we construct the main result of this section, which contains a quantitative characterization for strong minima. A similar result for group-sparsity minimization problem is recently established in \cite[Theorem~5.3]{FNT21}. 

\begin{Theorem}[Characterizations for strong minima of  low-rank optimization problems]\label{TheoNu} Suppose that $\OX\in \Phi^{-1}({\rm int}(\dom h))$ is a minimizer of  problem \eqref{Nucl} and $\OY=-\frac{1}{\mu}\Phi^*\nabla h(\Phi \OX)$ with decomposition \eqref{OX} and \eqref{OY2}. The following are equivalent:

\begin{enumerate}
    \item[{\rm(i)}] $\OX$
 is a strong solution to problem \eqref{Nucl}.
     \item[{\rm (ii)}] $\Ker \Phi\cap \mathcal{E}\cap \mathcal{C}=\{0\}$ with \begin{eqnarray}\label{C}
        &&\mathcal{E}\eqdef\left\{W\in \R^{n_1\times n_2}|\; P_{\TT}W\in \OU\begin{pmatrix}A& B&0\\B^T&0&0\\0&0&0\end{pmatrix}\OV^T,A\in \mathbb{S}^r, B\in \R^{r\times (p(\OY)-r)} \right\},  \label{E1}\\
    &&\mathcal{C}\eqdef\left\{W\in \R^{n_1\times n_2}|\; \la E,             W\ra+\|P_{\TT^\perp} W\|_*= 0\right\}\qquad \mbox{with}\qquad E\eqdef U_IV_I^T.   \label{C0}
    \end{eqnarray}
\item[{\rm (iii)}] The following  conditions (a) and either (b) or (c)  are satisfied: \begin{itemize}
         \item[{\rm (a)}] {\rm (Strong Restricted Injectivity)}: $\Ker \Phi\cap \mathcal{E}\cap \TT=\{0\}$.
         \item[{\rm (b)}] {\rm (Strong Nondegenerate Source Condition)}: There exists  $Y\in \Im \Phi^*+\mathcal{E}^\perp$ such that $Y=U \begin{pmatrix} \Id_r &0\\0 & Z\end{pmatrix} V^T $ and  $\|Z\|< 1$.
         \item [{\rm (c)}] {\rm (Analysis Strong Source Condition)}: The {\rm Strong Source Coefficient} $\zeta(\OX)$, which is the optimal value of the following spectral norm optimization problem
         \begin{equation}\label{Spec}
\min_{Z\in \R^{(m-r)\times (n-r)}}\qquad \|Z\|\qquad \mbox{subject to}\qquad \mathcal{M}(U_JZV_K^T)=-\mathcal{M}E
\end{equation}
is smaller than $1$, where $\mathcal{M}$ is a linear operator such that  $\Im \mathcal{M}^*=\Ker \Phi\cap \mathcal{E}$.

     \end{itemize}
\end{enumerate}
\end{Theorem}
\noindent{\bf Proof.} Let us verify the equivalence between (i) and (ii).  By Corollary~\ref{CoNu}, it suffices to show that 
\begin{equation}\label{Equi2}
    \Ker \Phi\cap T_{N_\B(\OY)}(\OX)=\Ker \Phi\cap \mathcal{E}\cap \mathcal{C}. 
\end{equation}
By Proposition \ref{Pro4} and Corollary~\ref{CoNu}, we have 
\begin{equation}\label{DecomT}
    T_{N_\B(\OY)}(\OX)=\mathcal{E}\cap \mathcal{C}(\OX,\OY). 
\end{equation}
As $\OY\in \Im \Phi^*$, note from \eqref{CritY} and \eqref{CXY} that 
\[
\Ker \Phi \cap \mathcal{C}(\OX,\OY)= \Ker \Phi\cap \mathcal{C}. 
\]
This together with  \eqref{DecomT}  verifies \eqref{Equi2}  and also the equivalence between (i) and (ii). 

Next, let us verify the implication [(ii)$\Rightarrow$(iii)].  Suppose that (ii) (or (i)) is satisfied. Note from the projection formula \eqref{Proj} that $\mathcal{E}\cap \TT\subset T_{N_\B(\OY)}(\OX)$. It follows from  Corollary~\ref{CoNu} that the Strong Resitricted Injectivity holds.

Since $\OY\in \partial \|\OX\|_*\cap \Im \Phi^*$, we obtain from \eqref{ConV} and \eqref{CXY} that 
\[
\la E,W\ra+\|P_{\TT^\perp} W\|_*=dg(\OX)(W)\ge \la \OY,W\ra= 0 \quad \mbox{for any}\quad W\in \Ker \Phi. 
\]
Condition (ii) means 
\[
c \eqdef \min\{\la E,W\ra+\|P_{\TT^\perp} W\|_*|\; W\in \Ker \Phi\cap \mathcal{E}, \|W\|_*=1\}>0,
\]
which implies that 
\begin{equation}\label{c}
    k(W)\eqdef\la E,W\ra+\|P_{\TT^\perp}W\|_*\ge c\|W\|_*\quad \mbox{for all}\quad W\in \Ker \Phi \cap \mathcal{E}.
\end{equation}

As $\Im \mathcal{M}^*=\Ker\Phi\cap \mathcal{E}$, for any $W\in \Ker\Phi\cap \mathcal{E}$, we write $W=\mathcal{M}^*Y$ and derive from \eqref{c} that 
\begin{equation}\label{MY}
(1-c)\|P_{\TT^\perp}\MM^*Y\|_*\ge\|P_{\TT^\perp}\MM^*Y\|_* -c\|\MM^*Y\|_*\ge -\la E, \MM^*Y\ra=-\la \MM E,Y\ra. 
\end{equation}
Since $\Im \mathcal{M}^*\subset \Ker\Phi$, we have $ \Im \Phi^*\subset \Ker \mathcal{M}$ and thus $\OY\in \Ker\mathcal{M}$. It follows that 
\[
0=\MM\OY=\MM P_{\TT}\OY+\MM P_{\TT^\perp}\OY=\MM E+\MM P_{\TT^\perp}\OY.
\]
This together with \eqref{MY} implies that 
\[
(1-c)\|P_{\TT^\perp}\MM^*Y\|_*\ge\la \MM P_{\TT^\perp}\OY,Y\ra=\la P_{\TT^\perp}\OY, \MM^*Y\ra=\la P_{\TT^\perp}\OY, P_{\TT^\perp}\MM^*Y\ra. 
\]
Define $\B_*\eqdef\{W\in \R^{m\times n}|\; \|W\|_*\le 1\}$ the unit ball with respect to nuclear norm. We obtain from the later and the classical minimax theorem \cite[Corollary 37.3.2]{R70} that 
\[\begin{array}{ll}
1-c&\disp\ge \sup_{W\in \B_*} \la P_{\TT^\perp}\OY, W\ra-\iota_{{\rm Im}\, P_{\TT^\perp}\MM^*}(W)\\
 &=\disp\sup_{W\in \B_*}\inf_{X\in {\rm Ker}\,  \MM P_{\TT^\perp}} \la P_{\TT^\perp}\OY+X, W\ra\\
 &=\disp\inf_{X\in {\rm Ker}\,  \MM P_{\TT^\perp}}\sup_{W\in \B_*} \la P_{\TT^\perp}\OY+X, W\ra\\
 &=\disp\inf_{X\in {\rm Ker}\,  \MM P_{\TT^\perp}}\| P_{\TT^\perp}\OY+X\|.
\end{array}
\]
Hence there exists $X_0\in \Ker \MM P_{\TT^\perp}$ such that $\| P_{\TT^\perp}\OY+X_0\|<1$. Due to the projection formula \eqref{Proj}, observe that
\[
1>\| P_{\TT^\perp}\OY+X_0\|\ge \| P_{\TT^\perp}(P_{\TT^\perp}\OY+X_0)\|=\|P_{\TT^\perp}(\OY+X_0)\|.
\]
Define $Y_0=\OY+P_{\TT^\perp}X_0$, we have
\[
\MM Y_0=\MM \OY+\MM P_{\TT^\perp}X_0=0.
\]
Note that  $\Ker \MM=\Im \Phi^*+\mathcal{E}^\perp$, $\OY\in \Im \Phi^*\subset \Ker \MM$, and $X_0\in \Ker \MM P_{\TT^\perp}$. It follows that $Y_0\in \Ker \MM=\Im \Phi^*+\mathcal{E}^\perp$.  Moreover, observe that $P_{\TT}Y_0=P_{\TT}\OY=E$ and  $\|P_{\TT^\perp}Y_0\|<1$. Thus, $Y_0$ satisfies the condition in (b). As $\Ker \MM=\Im \Phi^*+\mathcal{E}^\perp$,  (b) and (c) are equivalent, we ensure the implication [(ii)$\Rightarrow$(iii)].

It remains to justify the implication [(iii)$\Rightarrow$(ii)]. Suppose that the Strong Restricted Injectivity (a)  and the Strong Source Condition (b) hold with some $Y_0\in \Im \Phi^*+\mathcal{E}^\perp$ satisfying the condition in (b). Indeed, pick any $W\in \Ker \Phi\cap \mathcal{E}\cap\mathcal{C}=\Im^*\MM$. As $Y_0\in \Ker \MM$, we have $\la Y_0,W\ra=0$. It follows that  \[\begin{array}{ll}
0=\la E,W\ra+\|P_{\TT^\perp}W\|_*&=\disp\la P_{\TT}Y_0,W\ra+\|P_{\TT^\perp}W\|_*\\
&=\disp\la Y_0,W\ra-\la P_{\TT^\perp}Y_0,W\ra+\|P_{\TT^\perp}W\|_*\\
&=\disp -\la P_{\TT^\perp}Y_0,P_{\TT^\perp}W\ra+\|P_{\TT^\perp}W\|_*\\
&\ge  (1-\|P_{\TT^\perp}Y_0\|)\|P_{\TT^\perp}W\|_*. 
\end{array}
\]
Since $\|P_{\TT^\perp}Y_0\|<1$, we have $P_{\TT^\perp}W=0$, i.e., $W\in \TT$. This implies that $W=0$ due to the Strong Restricted Injectivity (a). The proof is complete.   \endproof

\begin{Remark}{\rm  The Strong Restricted  Injectivity (a) means that the linear operator $\Phi$ is injective on the subspace $\mathcal{E}\cap \TT$. It is similar with the one with the same name  in \cite[Theorem~5.3]{FNT21} that is used to characterize the uniqueness (strong) solution for {\em group-sparsity} optimization problems.   This condition is adopted from the Restricted Injectivity \eqref{RI1} in \cite{GHS11}; see also \cite{CP10,CR09,CR13} for the case of nuclear norm minimization problems.  The Strong Restricted  Injectivity is certainly 
weaker than the Restricted Injectivity. 

The Strong Nondegenerate Source Condition and Analysis Strong Source Condition also inherit the same terminologies introduced in \cite[Theorem~5.3]{FNT21} for group-sparsity optimization problems. In \cite{CP10,CR13,GHS11,VGFP15, VPF15}, the Nondegenerate Source Condition at $\OX$ means the existence of a {\em dual certificate} $Y_0\in \Im\Phi^*\cap\partial\|\OX\|_*$ satisfying  $\|P_{\TT^\perp}Y_0\|<1$, which is equivalent to
\begin{equation}\label{NSC}
    \Im \Phi^*\cap {\rm ri}\, \partial \|\OX\|_*\neq \emptyset. 
\end{equation}
This condition is weaker than the Nondegeneracy Condition \eqref{NC1}.   
In the case of $\ell_1$ optimization, it is well-known that the Restricted Injectivity and Nondegenerate Source Condition together characterize solution uniqueness \cite{GHS11}. For nuclear norm minimization problem, \cite{CR09,CR13} shown that they are  sufficient  for solution uniqueness of problems \eqref{Nu0} and \eqref{Nucl}; see also \cite{VGFP15,VPF15} for more general convex optimization problems. It is worth noting that they are not necessary conditions for solution uniqueness; see, e.g., \cite[Example~4.15]{FNT21}. One hidden reason is that the nuclear norm is not {\em polyhedral}. Recently  \cite{FNT21} shows that these two conditions characterize the so-called {\em sharp minima} of \eqref{NNM}, which is somewhat between solution uniqueness and strong minima; see our Remark~\ref{Sharp} for further discussion. Due to \cite[Proposition~4.27]{FNT21}, they are sufficient for strong solution of problem \eqref{Nucl}. This fact can be obtained from Theorem~\ref{TheoNu} by observing that our  Strong Nondegenerate Source Condition is weaker than Nondegenerate Source Condition, since  Strong Nondegenerate Source Condition involving the set $\mathcal{E}$  means
\begin{equation}\label{SNSC}
(\Im \Phi^*+\mathcal{E}^\perp)\cap {\rm ri}\, \partial \|\OX\|_*\neq\emptyset     
\end{equation} 
due to \eqref{Rei}.
}
\end{Remark}

\begin{Remark}[Checking  Strong Restricted Injectivity,  Strong Sufficient Condition \eqref{Suff}, and constructing the linear operator $\mathcal{M}$]\label{Comzeta} {\rm To check the Strong Restricted Injectivity,  observe first that 
\[
\mathcal{E}=\left\{\OU\begin{pmatrix}A& B&0\\B^T&C&D\\0&E&F\end{pmatrix}\OV^T\in \R^{n_1\times n_2}|\;A\in \mathbb{S}^r, B\in \R^{r\times (p-r)} \right\},
\]
which is a subspace of $\R^{n_1\times n_2}$ with dimension $q\eqdef\frac{r(r+1)}{2}+r(p-r)+(n_1-r)(n_2-r)$. Moreover, the restriction of $\mathcal{E}$ on $\TT$ 
\begin{equation}\label{ET}
\mathcal{E}\cap \TT=\left\{\OU\begin{pmatrix}A& B&0\\B^T&0&0\\0&0&0\end{pmatrix}\OV^T\in \R^{n_1\times n_2}|\;A\in \mathbb{S}^r, B\in \R^{r\times (p-r)} \right\},
\end{equation}
 is also a subspace of $\R^{n_1\times n_2}$ with dimension $s\eqdef\frac{r(r+1)}{2}+r(p-r)$. The set in $ \OU\begin{pmatrix}\mathbb{S}^{p}&0\\0&0\end{pmatrix}\OV^T $ in Strong Sufficient Condition \eqref{Suff} 
is another subspace of $\R^{m\times n}$ with dimension  $l\eqdef\frac{1}{2}p(p+1)$. 
Suppose that  $\{W_1, \ldots, W_s\}$ form a basis of $\mathcal{E}\cap \TT$,  $\{W_1, \ldots, W_l\}$ form a basis of $\OU\begin{pmatrix}\mathbb{S}^{p}&0\\0&0\end{pmatrix}\OV^T$ and $\{W_1\ldots W_{q}\}$ form a basis of $\mathcal{E}$. 
For any $W\in \mathcal{E}$, we write $W=\lm_1W_1+\ldots +\lm_qW_q$ and obtain that 
\begin{equation}\label{linear}
\Phi(W)=\lm_1\Phi(W_1)+\ldots +\lm_q\Phi(W_q).
\end{equation}
Define $\Psi\eqdef\begin{pmatrix}\Phi(W_1) &\ldots&\Phi(W_q)\end{pmatrix}$ to be an $m\times q$ matrix and $\Psi_s$, $\Psi_l$ to be the submatrices of the first $s,l$ columns of $\Psi$, respectively. By \eqref{linear}, the Strong Restricted Injectivity is equivalent to the condition that $\Ker \Psi_s=0$, i.e., $\rank{\Psi_s}=s.$ Similarly, the Strong Sufficient Condition \eqref{Suff} means $\rank{\Psi_l}=l$.

Next let us discuss how to construct the operator $\mathcal{M}$.  Let $\Hat U\Lm \Hat V^T$ be a SVD of $\Psi$ with $k={\rm rank}\, \Psi$. Define $\Hat V_G$ to be the $q\times (q-k)$ submatrix of $\Hat V$ where $G=\{k+1, \ldots,q\}$. Note that 
\[
\Im \Hat V_G=\Ker \Psi.
\]
Determine the linear operator $\mathcal{M}:\R^{n_1\times n_2}\to \R^{q-k}$ by
\begin{equation}\label{M}
    \mathcal{M}X\eqdef\Hat V_G^T\begin{pmatrix}\la W_1,X\ra,\ldots,\la W_q,X\ra\end{pmatrix}^T\quad\mbox{for all}\quad X\in \R^{n_1\times n_2}. 
\end{equation}
It is easy to check that $\Im \mathcal{M}^*=\Ker \Phi\cap\mathcal{E}$. 
   \endproof  }
\end{Remark}

\begin{Remark}[Checking the Analysis Strong Source Condition without using the linear operator $\mathcal{M}$]{\rm  To verify the Analysis Strong Source Condition, Theorem~\ref{TheoNu} suggests to solve the optimization problem \eqref{Spec}, which involves the linear operator $\mathcal{M}$. To avoid the computation of $\mathcal{M}$, note first that the constraint in \eqref{Spec} means $U_JZV_K^T+E\in \Ker \mathcal{M}=\Im \Phi^*+\mathcal{E}^\perp.$ Suppose that $\mathcal{N}$ is the linear operator with $\Ker \mathcal{N}=\Im \Phi^*$, i.e., $\Im \mathcal{N}^*=\Ker \Phi$, the later condition is equivalent to 
\[
\mathcal{N}(U_JZV_K^T+E+W)=0\quad \mbox{for some}\quad W\in \mathcal{E}^\perp,
\]
where $\mathcal{E}^\perp$ is computed by 
\begin{equation}\label{Eperp}
    \mathcal{E}^\perp=\left\{U\begin{pmatrix}A&B&C\\-B^T&0&0\\D&0&0\end{pmatrix}V^T|\;A\in \mathbb{V}_r, B\in \R^{r\times (p-r)}, C\in \R^{r\times (n_2-p)}, D\in \R^{ (n_1-p)\times r} \right\},
\end{equation}
which is a subspace of $\TT$ with dimension $\frac{1}{2}r(r-1)+r(m+n-p-r)$; here $\mathbb{V}_r\subset \R^{r\times r}$ is the set of all skew-symmetric matrices. Hence  problem \eqref{Spec} is equivalent to the following one 
\begin{equation}\label{Spec2}
\min_{Z\in \R^{(n_1-r)\times (n_2-r)}, W\in \R^{n_1\times n_2}}\qquad \|Z\|\quad \mbox{subject to}\quad \mathcal{N}(U_JZV_J^T+W)=-\mathcal{N}E\quad \mbox{and}\quad   W\in \mathcal{E}^\perp. 
\end{equation}
The price to pay is the size of this optimization problem is bigger than the one in \eqref{Spec}. But finding the linear operator $\mathcal{N}$ is much easier than $\mathcal{M}$. 
}
\end{Remark}

\begin{Corollary}[Quantitative sufficient condition for strong solution]\label{QSC} Suppose that  $\OX\in \Phi^{-1}({\rm int}(\dom h))$ is a minimizer of problem \eqref{Nucl}. Then $\OX$ is a strong solution provided that the Strong Restricted Injectivity holds and 
    \begin{equation}\label{GG}
    \gg(\OX)\eqdef \|P_{\TT^\perp}\MM^*(\MM P_{\TT^\perp}\MM^*)^{-1}\MM E\|<1,
\end{equation}
where $\mathcal{M}$ is the linear operator satisfying $\Im \mathcal{M}^*=\Ker\Phi\cap \mathcal{E}.$

\end{Corollary}
{\noindent \bf Proof.} Suppose that Strong Restricted Injectivity holds at the minimizer $\OX$. We consider the following linear equation:
\begin{equation}\label{LE}
    \MM P_{\TT^\perp} Y=-\MM E\quad \mbox{for}\quad Y\in \R^{n_1\times n_2}.
\end{equation}
As $\OY\in \partial\|\OX\|_*\cap \Im \Phi^*\subset \partial\|\OX\|_*\cap \Ker \MM$, we have 
\[
0=\MM\OY=\MM(P_{\TT}\OY+P_{\TT^\perp}\OY)=\MM E +P_{\TT^\perp}\OY.
\]
It follows that $\OY$ is a solution to \eqref{LE}. Another solution of \eqref{LE} is 
\begin{equation}\label{MPS}
    \Hat Y\eqdef-(\MM P_{\TT^\perp})^\dag\MM E, 
\end{equation}
where $(\MM P_{\TT^\perp})^\dag$ is the Moore-Penrose generalized inverse of $\MM P_{\TT^\perp}$. Next we claim that $\MM^*P_{\TT^\perp}\MM^*:\Im \MM\to \Im \MM$ is a bijective mapping. Indeed, suppose that $\MM P_{\TT^\perp}\MM^*z=0$ for some $z\in \Im \MM$, we have 
\[
\|P_{\TT^\perp}\MM^*z\|^2=\la P_{\TT^\perp}\MM^*z,P_{\TT^\perp}\MM^*z\ra =0.
\]
It follows that  $M^*z\in \TT\cap\Im\MM^*=\{0\}$, which implies $z\in \Ker\MM^*\cap\Im \MM=\{0\}$. Thus  $\MM P_{\TT^\perp}\MM^*$ is injective in $\Im \MM$. As the operator is self-dual, it is bijective. We obtain that  
\[
\Hat Y=-P_{\TT^\perp}\MM^* (\MM P_{\TT^\perp} \MM^*)^{-1}\MM E\in \TT^\perp.
\]
By the decomposition \eqref{Proj}, we may write
\[
\Hat Y=U\begin{pmatrix} 0&0\\0&\Hat Z\end{pmatrix}V^T\quad \mbox{for some}\quad \Hat Z\in \R^{(n_1-r)\times (n_2-r)}.
\]
Observe from \eqref{Spec} that 
\[
\zeta(\OX)\le \|\Hat Z\|=\|\Hat Y\|=\gg(\OX). 
\]
If $\gg(\OX)<1$, we  have $\zeta(\OX)<1.$ $\OX$ is a strong solution due to Theorem~\ref{TheoNu}.\endproof

\section{Characterizations for strong minima of  nuclear norm minimization problems}
\setcounter{equation}{0}

In this section, let us consider the nuclear norm minimization problem \eqref{Nu0}:
\begin{equation}\label{NNM}
    \min_{X\in \R^{n_1\times n_2}}\quad \|X\|_*\quad \mbox{subject to}\quad \Phi X=M_0,
\end{equation}
where $\Phi:\R^{n_1\times n_2}\to \R^m$ is a linear operator ($n_1\le n_2$), $M_0\in \R^m$ is a known observation. This is a particular case of problem \eqref{CO} with $g(X)=\|X\|_*$ for $X\in \R^{n_1\times n_2}$. Note that $X_0$ is a solution of this problem if and only if $\Phi X_0=M_0$ and
\[
0\in \partial \|X_0\|_*+N_{\Phi^{-1}(M_0)}(X_0)=\partial \|X_0\|_*+\Im \Phi^*,
\]
which is equivalent to the existence of a {\em dual certificate} $Y\in   \Im \Phi^*\cap \partial\|X_0\|_*$. Define 
\[
\Delta(X_0)\eqdef \Im \Phi^*\cap \partial\|X_0\|_*\]
to be the set of all dual certificates.  
\begin{Lemma}\label{Iso} Let $\Omega$ be a nonempty closed set of $\R^{n_1\times n_2}$. Suppose that $U_1,U_2\in \R^{n_1\times n_2}$ and $V_1,V_2\in \R^{n_2\times n_2}$ are orthogonal matrices satisfying $U_1\Omega V_1^T\supset U_2\Omega V_2^T$. Then we have 
\begin{equation}\label{Om}
U_1\Omega V_1^T= U_2\Omega V_2^T.
\end{equation}
\end{Lemma}
{\noindent Proof.} Define $U=U_1^TU_2$ and $V=V_1^TV_2$. It is easy to check that $U$ and $V$ are orthogonal matrices. We have 
\begin{equation}\label{sub}
 U\Omega V^T\subset \Omega.
\end{equation}
Define the mapping $\ph:\Omega\to\Omega$ by $\ph(X)=UXV^T$ for all $X\in \Omega$. Note that $\ph$ is an {\em isometry} with the spectral metric $d(X,Y)=\|X-Y\|$ for all $X,Y\in \Omega$. Indeed, we have
\[
d(\ph(X),\ph(Y))=\|U(X-Y)V^T\|=\|X-Y\|\quad \mbox{for all}\quad X,Y\in \Omega. 
\]
Note further that $\|\ph(X)\|=\|X\|$ for all $X\in \Omega$. It follows from \eqref{sub} that 
\[
\ph(\Omega\cap \Bar\B_k(0))\subset \Omega\cap \Bar \B_k(0)\quad \mbox{for any}\quad k\in \mathbb{N}. 
\]
Since $\Omega\cap \Bar \B_k(0)$ is a compact metric space and $\ph$ is an isometry, it is well-known that 
\[
\ph(\Omega\cap \Bar \B_k(0))= \Omega\cap \Bar \B_k(0).
\]
Hence we have
\[
\ph(\Omega)= \ph\left(\bigcup_{k=1}^\infty  (\Omega\cap \Bar \B_k(0))\right)= \bigcup_{k=1}^\infty \ph\left( \Omega\cap \Bar \B_k(0)\right)=\bigcup_{k=1}^\infty  (\Omega\cap \Bar \B_k(0))=\Omega. 
\]
This verifies \eqref{Om} and completes the proof of the lemma. \endproof

For any $Y\in \partial \|X_0\|_*$, recall from  \eqref{p} that $p(Y)$ to be the number of singular values of $Y$ that are equal to $1$ . It follows from Lemma~\ref{MLem} that $r\eqdef\rank(X_0)\le p(Y)\le n_1$. Define the following constant 
\begin{equation}\label{q}
q(X_0)\eqdef\min\{p(Y)|\; Y\in \Delta(X_0)\}.
\end{equation}
When $X_0$ is a minimizer of problem \eqref{NNM}, $q(X_0)$ is well-defined and bigger than or equal to $\rank(X_0)$. The following theorem is the main result in this section. 

\begin{Theorem}[Characterizations for strong minima of nuclear norm minimization problems]\label{SSN} Suppose that $X_0$ is an optimal solution of problem \eqref{NNM}. The following are equivalent:

{\bf (i)} $X_0$ is a strong solution of  problem \eqref{NNM}.

{\bf (ii)}  The following equality holds
\begin{equation}\label{ma}
    \bigcap_{Y\in \Delta(X_0)}\Big[\Ker \Phi\cap  T_{N_\B(Y)}(X_0)\Big]=\{0\}.
\end{equation}

{\bf (iii)} There exists a dual certificate  $\OY\in \Delta(X_0)$ such that 
\begin{equation}\label{ma2}
    \Ker \Phi\cap  T_{N_\B(\OY)}(X_0)=\{0\}.
\end{equation}

{\bf (iv)} For any   $Y\in \Delta(X_0)$ satisfying $p(Y)=q(X_0)$, condition \eqref{ma2} is satisfied.
\end{Theorem}
{\noindent \bf Proof.} Note that the nuclear norm function $\|\cdot\|_*$ is second order regular at $X_0$ \cite{CDZ17} and it also satisfies the quadratic grown condition at $X_0$, as $\partial \|\cdot\|_*$ is metrically subregular at $X_0$ for any $Y\in \partial \|X_0\|_*$; see \cite{ZS17}. By applying  Lemma~\ref{Lem} and Theorem~\ref{thm5} with $\XX=\R^{n_1\times n_2}$, $\YY=\R^{m}$, $g(\cdot)=\|\cdot\|_*$, and $K=\{M_0\}$, we have $X_0$ is a strong solution if and only if 
\begin{equation}\label{Equi3}
    \left[\bigcap_{Y\in \Delta(X_0)}T_{N_\B(Y)}(X_0)\right]\cap C(X_0)=\{0\},
\end{equation}
where $C(X_0)$ is the critical cone \eqref{Cr} computed by
\[
C(X_0)=\{W\in \R^{n_1\times n_2}|\; W\in \Ker \Phi, dg(X_0)(W)=0\}. 
\]
By Proposition~\ref{Pro4} and formula~\eqref{Tang2}, we note that if  $W\in \Ker \Phi\cap T_{N_\B(Y)}(X_0)$ for some $Y\in \Delta(X_0)$ then $W\in \mathcal{C}(X_0,Y)$, i.e., $dg(X_0)(W)=\la Y,W\ra=0$. Thus $W\in C(X_0)$ and that 
\[
\left[\bigcap_{Y\in \Delta(X_0)}T_{N_\B(Y)}(X_0)\right]\cap C(X_0)=\left[\bigcap_{Y\in \Delta(X_0)}T_{N_\B(Y)}(X_0)\right]\cap \Ker \Phi.
\]
This together with \eqref{Equi3} verifies the equivalence between (i) and (ii). 

The implication [(iii)$\Rightarrow$(ii)] and [(iv)$\Rightarrow$(ii)] are trivial. To justify the converse implications, we first claim the existence of some dual certificate $\OY\in \Delta(X_0)$  such that  $p(\OY)=q(X_0)$ and 
\begin{equation}\label{Zorn}
\bigcap_{Y\in \Delta(X_0)}\Big[  T_{N_\B(Y)}(X_0)\Big]=T_{N_\B(\OY)}(X_0).
\end{equation}
We prove this by using a popular version  Zorn's lemma \cite[Proposition~5.9]{C12} on {\em partially directed order set}. Consider  the following  partially ordered set ({\em poset})
\begin{equation}\label{PO}
\mathcal{P}=\Big\{T_{N_\B(Y)}(X_0)|\; Y\in \Delta(X_0)\Big\} 
\end{equation}
with the  partial ordering  $\supset$ between sets in $\mathcal{P}$. Take into account any {\em downward chain}
\begin{equation}\label{chain}
T_{N_\B(Y_1)}(X_0)\supset T_{N_\B(Y_2)}(X_0)\supset \ldots \supset T_{N_\B(Y_k)}(X_0)\supset\ldots
\end{equation}
for a sequence $\{Y_k\}\subset \Delta(X_0)$. We can find a subsequence $\{Y_{k_l}\}$ of $\{Y_k\}$ such that they have the same rank $p\ge r$. According to the computation \eqref{Tang2}, there exist  orthogonal matrices $(U_{k_l},V_{k_l})\in\mathcal{O}(Y_{k_l})\cap \mathcal{O}(X_0)$ such that 
\[
T_{N_\B(Y_{k_l})}(X_0)=U_{k_l}\Omega_p V^T_{k_l}\quad\mbox{with}\quad  \Omega_p\eqdef\left\{\begin{pmatrix}A&B&0\\B^T&C&0\\0&0&0\end{pmatrix}|\; A\in \mathbb{S}^r, B\in \R^{r\times (p-r)}, C\in \mathbb{S}_+^{p-r}\right\}.
\]
It follows that 
\begin{equation}\label{sup2}
T_{N_\B(Y_{k_l})}(X_0)=U_{k_l}\Omega_p V^T_{k_l}\supset U_{k_{l+1}}\Omega_p V^T_{k_{l+1}}=T_{N_\B(Y_{k_{l+1}})}(X_0).
\end{equation}
Since $\Omega_p$ is closed, we obtain from Lemma~\ref{Iso} that 
\[
T_{N_\B(Y_{k_l})}(X_0)=T_{N_\B(Y_{k_{l+1}})}(X_0)\quad \mbox{for all}\quad l=1,2,\ldots.
\]
Hence the chain \eqref{chain} is bounded below by $T_{N_\B(Y_{k_1})}(X_0)\in \mathcal{P}$. This means that every downward chain of of $\mathcal{P}$ has a minimum in $\mathcal{P}$.

Let us show next that the poset $\mathcal{P}$ is  {\em directed downward} with the partial ordering  $\supset$ in the sense that for any two elements  $T_{N_\B(Y_1)}(X_0)$ and $T_{N_\B(Y_2)}(X_0)$ of $\mathcal{P}$ with $Y_1,Y_2\in \Delta(X_0)$, there exists $Y_3\in \Delta(X_0)$ such that 
\begin{equation}\label{RM}
T_{N_\B(Y_1)}(X_0)\supset T_{N_\B(Y_3)}(X_0)\quad\mbox{and}\quad T_{N_\B(Y_2)}(X_0)\supset T_{N_\B(Y_3)}(X_0).
\end{equation}
Indeed, we choose $Y_3=\frac{1}{2}(Y_1+Y_2)\in \Delta(X_0)$. For any $W\in T_{N_\B(Y_3)}(X_0)$,  we obtain from  \eqref{Tang2} that $0=d^2g(X_0|Y_3)(W)$. Hence, there exists $t_k\dn 0$ and $W_k\to W$ such that 
\[\begin{array}{ll}
 \dfrac{1}{k}&\ge \dfrac{g(X_0+t_kW_k)-g(X_0)-0.5\la Y_1+Y_2,X_0\ra}{0.5t_k^2}\\
 &=\dfrac{g(X_0+t_kW_k)-g(X_0)-\la Y_1,X_0\ra}{t_k^2}+\dfrac{g(X_0+t_kW_k)-g(X_0)-\la Y_2,X_0\ra}{t_k^2},
 \end{array}
\]
which implies that \[
0\ge d^2g(X_0|Y_1)(W)+d^2g(X_0|Y_2)(W).\]
Since   $d^2g(X_0|Y_1)(W), d^2g(X_0|Y_2)(W)\ge 0$, we obtain from Lemma~\ref{Dom} and Lemma~\ref{Lem} that 
\[
W\in \Ker d^2g(X_0|Y_1)=T_{N_\B(Y_1)}(X_0)\quad \mbox{and}\quad W\in \Ker d^2g(X_0|Y_2)=T_{N_\B(Y_2)}(X_0).
\]
This clearly verifies the directed condition \eqref{RM}. By \cite[Proposition~5.9]{C12}, the directed downward poset $\mathcal{P}$ has a minimum in $\mathcal{P}$ in the sense that there exists $\OY$ such that \eqref{Zorn} is valid. The implication [(ii)$\Rightarrow$(iii)] follows from \eqref{Zorn}. 

Next, let us show that $p(\OY)=q(X_0)$ and 
\begin{equation}\label{Zorn2}
T_{N_\B(Y)}(X_0)= T_{N_\B(\OY)}(X_0)
\end{equation}
for any $Y\in \Delta(X_0)$ with $p(Y)=q(X_0)$. Indeed, pick any $Y$ satisfying the later condition, we obtain from \eqref{Zorn} that
\begin{equation}\label{Sup3}
T_{N_\B(Y)}(X_0)\supset T_{N_\B(\OY)}(X_0).
\end{equation}
Due to the computation \eqref{Tang2}, the maximum ranks of matrices in $T_{N_\B(\OY)}(X_0)$ and $T_{N_\B(Y)}(X_0)$ are  $p(\OY)$ and $p(Y)$, respectively. It follows that $q(X_0)\le p(\OY)\le p(Y)= q(X_0)$, which verifies that $p(\OY)=q(X_0)$.
Moreover, due to \eqref{Tang2} and \eqref{Sup3}, we get \eqref{Zorn2} from Lemma~\ref{Iso}. Combining \eqref{Zorn2} and  \eqref{Zorn} ensures the implication  [(ii)$\Rightarrow$(iv)]. The proof  is complete. \endproof

\begin{Remark}{\rm A sufficient condition for strong minima of nuclear norm minimization \eqref{NNM} can be obtained from \cite[Theorem~12]{CDZ17}. However, their condition has a format of a minimax problem: for any element in the critical cone, there exists some Lagrange multiplier such that a certain second order sufficient condition is satisfied, i.e., the Lagrange multiplier used in the sufficient condition depends on the choice of elements in critical cone. This is a typical situation; see \eqref{MinMax} for instance.  In our characterizations for strong minima in part (iii) and (iv) of the above theorem, the existence of dual certificate $\OY$ 
is independent from elements of the critical cone. Condition (iii) is close to the {\em maximin} situation. Morever, we know extra information about these dual certificates from (iv) that they should have minimum  numbers of singular values that are equal to 1.

}\end{Remark}

Similarly to Theorem~\ref{TheoNu}, condition \eqref{ma2} is equivalent to the combination of Strong Restricted Injectivity and Strong Nondegenerate Source Condition. Let us recall   the model tangent subspace \eqref{T}  at $X_0$ here:
\[
\TT_0\eqdef\{U_0Y^T+XV_0^T|\; X\in \R^{n_1\times r}, Y\in \R^{n_2\times r}\},
\]
where $U_0\Sigma_0V_0^T$ is a compact SVD of $X_0$. 

\begin{Corollary}[Strong Restricted Injectivity and Strong Nondegenerate Source Condition for strong minima ]\label{CoroS} Suppose that $X_0$ is a minimizer of problem \eqref{NNM}. Then $X_0$ is a strong solution of problem \eqref{NNM} if and only if both following conditions hold:
\begin{itemize}
\item[(a)] {\rm Strong Restricted Injectivity:} $\Ker \Phi \cap \mathcal{E}_0\cap T_0=\{0\}$, where 
\begin{equation}\label{E_0}
\mathcal{E}_0=\left\{W\in \R^{n_1\times n_2}|\; P_{T_0}W\in U\begin{pmatrix}A& B&0\\B^T&0&0\\0&0&0\end{pmatrix}V^T,A\in \mathbb{S}^r, B\in \R^{r\times (q(X_0)-r)} \right\}
\end{equation}
with some $(U,V)\in \mathcal{O}(X_0)\cap\mathcal{O}(Y_0)$ with $Y_0\in \Delta(X_0)$ satisfying  $p(Y_0)=q(X_0)$. 

\item[(b)] {\rm Strong Nondegenerate Source Condition:} There exists $Y\in \Im\Phi^*+\mathcal{E}_0^\perp$ such that $Y\in {\rm ri}\, \partial \|X_0\|_*. $
\end{itemize}
\end{Corollary}

\begin{Remark}[Sharp minima vs Strong minima]\label{Sharp2}{\rm  
The set $\mathcal{E}_0$ is only dependent on $X_0$. Indeed, it follows from \eqref{Zorn} and \eqref{Tang2} that 
\[
\mathcal{E}_0=\left\{W\in \R^{n_1\times n_2}|\; P_{\TT}W\in P_{\TT}\left (\bigcap_{Y\in \Delta(X_0)}T_{N_\B(Y)}(X_0)\right)\right\},
\]
which is a subspace of $\R^{n_1\times n_2}$. As discussed after Theorem~\ref{TheoNu}, the Strong Restricted Injectivity is weaker than the Restricted Injectivity 
\begin{equation}\label{RI}
\Ker \Phi\cap \TT_0=\{0\}.   
\end{equation}
The Strong Nondegenerate Source Condition is also weaker than the  Nondegenerate Source Condition:
\begin{equation}\label{NSC2}
\Im \Phi^*\cap {\rm ri}\, \partial \|X_0\|_*\neq\emptyset.
\end{equation}
This condition means $q(X_0)=\rank (X_0)$. The Nondegenerate Source Condition together with  the Restricted Injectivity is used  in \cite{CR09,CR13} as sufficient conditions for solution uniqueness of problem \eqref{NNM}  at $X_0$. These two conditions are shown recently to be equivalent to the stronger property,  {\em sharp minima} at $X_0$  in \cite[Theorem~4.6]{FNT21} in the sense that there exists some $c>0$ such that 
\begin{equation}\label{Sharp}
\|X\|_*-\|X_0\|_*\ge c\|X-X_0\|\quad \mbox{for any $X\in \R^{n_1\times n_2}$ satisfying}\quad \Phi X=M_0. 
\end{equation}
Our Strong Restricted Injectivity and Strong Nongenerate Source Condition are   characterizations for a weaker property, the strong minima for problem \eqref{NNM}. Of courses, they can also serve as sufficient conditions for solution uniqueness at $X_0$; see \cite{HP23} for some recent  characterizations for this property.

In order to check Nondegenerate Source Condition \eqref{NSC2}, one has to show that the {\em Source Coefficient}  $\rho(X_0)$, the optimal value of the following optimization problem  
\begin{equation}\label{SC22}
    \min_{Z\in \TT_0^\perp}\quad \|Z\|\quad \mbox{subject to}\quad \mathcal{N}Z=-\mathcal{N}E_0\quad \mbox{with}\quad E_0\eqdef U_0V_0^T
\end{equation}
to be smaller than $1$, where $\mathcal{N}$ is a linear operator satisfying $\Ker \mathcal{N}=\Im \Phi^*$; see, e.g. \cite[Remark 4.5]{FNT21}. When the Restricted Injectivity \eqref{RI} holds, an upper bound for $\rho(X_0)$ is 
\begin{equation}\label{tau}
    \tau(X_0)\eqdef\|\mathcal{N}^*_{\TT_0^\perp}(\mathcal{N}^*_{\TT_0^\perp}\mathcal{N}_{\TT_0^\perp})^{-1} \mathcal{N}E_0\|\quad \mbox{with}\quad \mathcal{N}_{\TT_0^\perp}\eqdef \mathcal{N}P_{\TT_0^\perp}. 
\end{equation}
Hence condition $\tau(X_0)<1$ is sufficient for   sharp minima; see, e.g., \cite[Corollary~4.8]{FNT21}. This condition is known as the {\em Analysis Exact Recovery Condition} in \cite{NDEG13} for the case of $\ell_1$ optimization.  Another independent condition is also used to check solution uniqueness of  problem \eqref{NNM} is the so-called  {\em Irrepresentability Criterion} \cite{CR09,CR13,VPF17}:
\begin{equation}\label{IC}
{\bf IC}(X_0)\eqdef \|\Phi^*_{\TT_0^\perp}\Phi_{\TT_0}\left(\Phi_{\TT_0}^*\Phi_{\TT_0}\right)^{-1}E_0\|<1\quad \mbox{with}\quad \Phi_{\TT_0}\eqdef\Phi P_{\TT_0}.
\end{equation}
Note that ${\bf IC}(X_0)\ge\rho(X_0)$. Thus ${\bf IC}(X_0)<1$ also implies that $X_0$ is a sharp solution of problem \eqref{NNM}.
}
\end{Remark}

\begin{Remark}[Descent cone vs  tangent cone]\label{DC}{\rm Sharp minima and strong minima  of  problem \eqref{NNM} are sufficient for solution uniqueness, which is a significant property for exact recovery \cite{ALMT14,CRPW12,CR13}. An important geometric structure used to study solution uniqueness is the {\em descent cone} \cite{CRPW12} at $X_0$ defined by 
\begin{equation}\label{Des}
\mathcal{D}(X_0)\eqdef{\rm cone}\,\{X-X_0|\; \|X\|_*\le \|X_0\|_*\}.
\end{equation}
Indeed, \cite{CRPW12} shows that $X_0$ is a unique solution of problem \eqref{NNM} if and only if 
\begin{equation}\label{C12}
\Ker \Phi \cap \mathcal{D}(X_0)=\{0\}. 
\end{equation}
Unlike the tangent cones in \eqref{ma} or \eqref{ma2}, the descent cone $\mathcal{D}(X_0)$ may be not closed. Although the direct connection between the descent cone ${\mathcal D}(X_0)$ and tangent cones $T_{N_\B(Y)}(X_0)$ is not clear,  we claim  that 
\begin{equation}\label{DT}
\Ker \Phi \cap \mathcal{D}(X_0)\subset \bigcap_{Y\in \Delta(X_0)}\left[\Ker \Phi\cap T_{N_\B(Y)}(X_0)\right] 
\end{equation}
when $X_0$ is a minimizer of problem \eqref{NNM}, where $\Delta(X_0)=\Im  \Phi^*\cap \partial \|X_0\|_*$ is the set of dual certificates at $X_0$.
Indeed, for any $W\in \Ker \Phi \cap \mathcal{D}(X_0)$, there exists $\tau>0$ such that $\|X_0+\tau W\|_*\le \|X_0\|_*$. As $\Phi(X_0+\tau W)=\Phi X_0$, we have $\|X_0+\tau W\|_*= \|X_0\|_*$. Pick any $Y \in \Delta(X_0)$. In view of \eqref{Fen} and the definition of $Y$, it follows that 
\[
\|X_0+\tau W\|_*=\|X_0\|_*=\la Y, X_0\ra=\la Y, X_0+\tau W\ra,
\]
which implies that $X_0+\tau W\in N_\B(Y)$ due to \eqref{Fen}-\eqref{Inver} and thus $W\in T_{N_\B(Y)}(X_0)$. This verifies inclusion \eqref{DT}. Inclusion \eqref{DT} also tells us that  condition \eqref{ma} is sufficient for \eqref{C12}. This observation is not a surprise in the sense of Theorem~\ref{SSN}, as strong minima obviously implies solution uniqueness.  But solution uniqueness of problem \eqref{NNM} does not indicate strong minima; see \cite[Example~5.11]{FNT21}. Hence, inclusion \eqref{DT} is strict.   
}
\end{Remark}

Similarly to Corollary~\ref{Strict1}, the following result reveals the role of Strict Restricted Injectivity in strong minima. 

\begin{Corollary}[Strict Restricted Injectivity for strong minima of problem \eqref{NNM}]\label{Ri} Suppose that  $U_0\Sigma_0V_0^T$ is a compact SVD of  $X_0$ with $r={\rm rank}\,(X_0)$. If  $X_0$ is a strong solution, then the following  {\em Strict Restricted Injectivity} is satisfied: 
\begin{equation}\label{SS2}
    \Ker \Phi \cap U_0\mathbb{S}^rV_0^T=\{0\}.
\end{equation}

This condition is also  sufficient   for strong solution at $X_0$ provided that  {\em Nondegenerate Source Condition} \eqref{NSC2} holds at $X_0$. 
\end{Corollary}
\noindent{\bf Proof.} Note from \eqref{Tang2},
\[
T_{N_\B(Y)}(X_0)\supset U_0\mathbb{S}^rV_0^T \quad \mbox{for any}\quad Y\in \Delta(X_0).  
\]
If $X_0$ is a strong solution, combing \eqref{ma} with the latter inclusion verifies \eqref{SS2}. 

If Nongegenerate Source Condition \eqref{NSC2} holds at $X_0$, there exists $Y_0\in \Im \Phi^*\cap {\rm ri}\, \partial \|X_0\|_*$. Hence $\Delta(X_0)\neq \emptyset$, i.e., $X_0$ is a solution of problem \eqref{NNM}. It follows from Lemma~\ref{Lem} that  $p(Y_0)=\rank(X_0)$ and from \eqref{Tang2} that
\[
T_{N_\B(Y)}(X_0)=U_0\mathbb{S}^rV_0^T.
\]
Hence the validity of \eqref{SS2} implies that $X_0$ is a strong solution to problem \eqref{NNM} due to the equivalence between (i) and (iii) in Theorem~\ref{SSN}.  \endproof

Suppose that $U\begin{pmatrix}\Sigma_0&0\\0&0\end{pmatrix}V^T$ is a full SVD of $X_0$. The model tangent space $\TT_0$ at $X_0$ can be represented by 
\begin{equation}\label{T0}
\TT_0=\left\{U\begin{pmatrix}A&B\\C&0\end{pmatrix}V^T|\; A\in \R^{r\times r}, B\in \R^{r\times (n_2-r)}, C\in \R^{(n_1-r)\times r}\right\}.
\end{equation} 
It has dimension $r(n_1+n_2-r)$. When Restricted Injectivity \eqref{RI} holds, we have $m\ge r(n_1+n_2-r)$. Similarly, as the dimension of $U_0\mathbb{S}^rV_0^T$ is $\frac{1}{2}r(r+1)$, it is necessary for Strict Restricted Injectivity  \eqref{SS2} that $m\ge \frac{1}{2}r(r+1)$. Next we show that this bound $\frac{1}{2}r(r+1)$ for $m$ is tight.

\begin{Corollary}[Minimum bound for strong exact recovery] Suppose that $X_0$ is an $n_1\times n_2$ matrix with rank $r$. Then one needs at least $\frac{1}{2}r(r+1)$ measurements $m$ of $M_0$ so that solving the nuclear norm minimization problem \eqref{NNM} recovers exactly the strong solution $X_0$.

Moreover, there exist infinitely many linear operators $\Phi:\R^{n_1\times n_2}\to \R^{\frac{1}{2}r(r+1)}$ such that $X_0$ is a strong solution of problem \eqref{NNM}. 
\end{Corollary}
\noindent{\bf Proof.} Suppose that $U_0\Sigma V_0^T$ is a compact SVD of $X_0$. Let  $\{A_1, \ldots,A_s\}$ with $s\eqdef\frac{1}{2}r(r+1)$ be any basis of $U_0\mathbb{S}^rV_0^T$. If $X_0$ is a strong solution of problem \eqref{NNM}, Strict Restricted Injectivity \eqref{SS2} holds by Corollary~\ref{Ri}. It follows that $\{\Phi(A_1), \ldots,\Phi(A_s)\}$ are linearly independent. Hence, we have $m\ge s$, which verifies the first part. 

To justify the second part, we construct the linear operator $\Phi_s:\R^{n_1\times n_2}\to \R^s$ as follows:
\begin{equation}\label{Phi}
\Phi_s(X)\eqdef(\la A_k,X\ra)_{1\le k\le s}\in \R^s\quad \mbox{for any}\quad X\in \R^{n_1\times n_2}.
\end{equation}
Note that $\Im \Phi_s^*=\span\{A_1,\ldots,A_s\}=U_0\mathbb{S}^rV_0^T$. It follows that $E_0=U_0V_0^T\in \Im\Phi_s^*\cap \ri \partial \|X_0\|_*$ is a dual certificate that satisfies Nondegenerate Source Coundition \eqref{NSC2}. As 
\[
\Ker \Phi_s =(\Im \Phi_s^*)^\perp=(U_0\mathbb{S}^rV_0^T)^\perp,
\] we have $\Ker \Phi_s\cap U_0\mathbb{S}^rV_0^T=\{0\}$. Hence Strict Restricted Injectivity \eqref{SS2} holds. By Corollary~\ref{Ri} again, $X_0$ is a strong solution of problem \eqref{NNM}. \endproof 

\begin{Remark}[Low bounds on the number of measurement for exact recovery]{\rm In the theory of exact recovery, \cite{CRPW12} shows that with $m\ge 3r(n_1+n_2-r)$ random Gaussian measurements 
 it is sufficient to recover exactly $X_0$ with high probability by solving problem \eqref{NNM} from  observations $M_0=\Phi X_0$; see also \cite{CR13} for a similar result with a different approach. Also in \cite[Propositions~4.5 and 4.6]{CRPW12}, a lower bound on the number of measurements is discussed for exact recovery in {\em atomic norm} minimization via the descent cone \eqref{DC} and Terracini's Lemma \cite{H92}. The latter is used to get an estimate of the dimension of a subspace component of the descent cone. In the case of nuclear norm, this lower bound is indeed $\min\{n_1n_2, (r+1)(n_1+n_2)-r\}$; see aslo \cite[Proposition~12.2]{H92}. This bound holds for any linear measurement scheme. Our lower bound $\frac{1}{2} r(r+1)$ for $m$ is much smaller and only depends on the rank of $X_0$, but is achieved for special $\Phi$ only.  
}
\end{Remark}

\begin{Example}{\rm  Let us consider the following nuclear norm minimization problem
\begin{equation}\label{MC1}
    \min_{X\in \R^{2\times 2}}\|X\|_*\quad\mbox{subject to}\quad \Phi X\eqdef\begin{pmatrix}X_{11}& 0\\0&X_{22}\end{pmatrix}=\begin{pmatrix}1&0\\0&0\end{pmatrix},
\end{equation}
which is a matrix completion problem. Set $X_0=\begin{pmatrix}1&0\\0&0\end{pmatrix}$, we have 
\begin{equation}\label{PX}
\partial \|X_0\|_*=\begin{pmatrix}1&0\\0&[-1,1]\end{pmatrix}\quad \mbox{ and }\quad \TT_0=\left\{\begin{pmatrix}a&b\\c&0\end{pmatrix}|\; a,b,c\in \R\right\}.
\end{equation}
Moreover, note that 
\begin{equation}\label{KI}
\Ker \Phi=\left\{\begin{pmatrix}0&b\\c&0\end{pmatrix}|\; b,c\in \R\right\}\quad \mbox{ and }\quad \Im \Phi^*=\left\{\begin{pmatrix}a&0\\0&d\end{pmatrix}|\; a,d\in \R\right\}
\end{equation}
Note that $\Delta(X_0)=\Im\Phi^*\cap \partial\|X_0\|_*=\partial\|X_0\|_*$. Hence $X_0$ is a solution of problem \eqref{MC1}. However, $\Ker \Phi\cap \TT_0=\Ker \Phi\neq \{0\}$, i.e.,  Restricted Injectivity \eqref{RI} fails. Thus $X_0$ is not a sharp solution of problem \eqref{MC1}. Note that $q(X_0)=1$ and $Y_0=\begin{pmatrix}1&0\\0&0\end{pmatrix}\in \Delta(X_0)$ with $p(Y_0)=1$. We have from \eqref{Tang2} that 
\[
T_{N_{\B}(Y_0)}(X_0)=\left\{\begin{pmatrix}a&0\\0&d\end{pmatrix}|\; a,d\in \R\right\}.
\]
It is clear that $\Ker \Phi\cap T_{N_{\B}(Y_0)}(X_0)=\{0\}$. This shows that $X_0$ is a strong solution by Theorem~\ref{SSN}. \endproof
}
\end{Example}

\begin{Example}[Checking strong minima numerically]{\rm 
Let us consider the following matrix completion problem 
\begin{equation}\label{MC11}
\min_{X\in \R^{3\times 3}}\quad \|X\|_*\quad \mbox{subject to}\quad {\rm P}_\Omega(X)=M_0\eqdef\begin{pmatrix}4 &2& 4\\2&1&0\\4&0&0\end{pmatrix},
\end{equation}
where ${\rm P}_\Omega$ is the projection mapping defined by 
\[
{\rm P}_\Omega(X)\eqdef\begin{pmatrix}X_{11} &X_{12}& X_{13}\\X_{21}&X_{22}&0\\X_{31}&0&0\end{pmatrix}. 
\]
Define 
\[
X_0\eqdef \begin{pmatrix}4 &2& 4\\2&1&2\\4&2&4\end{pmatrix}\quad \mbox{and} \quad U=V\eqdef\frac{1}{3}\begin{pmatrix}2 &-2& 1\\1&2&2\\2&1&-2\end{pmatrix}.
\]
Note that $U,V$ are orthogonal matrices, ${\rm P}_\Omega(X_0)=M_0$,  $X_0=U \begin{pmatrix}9 &0& 0\\0&0&0\\0&0&0\end{pmatrix}V^T$ is an SVD of $X_0$, and ${\rm rank}\,(X_0)=1$. To check whether $X_0$ is a sharp solution of problem \eqref{MC11}, we compute the constants $\tau(E_0)$ in \eqref{tau} or $\rho(E_0)$ from problem \eqref{SC22}. In this case the linear operator $\mathcal{N}$ in \eqref{SC22} is chosen by $P_{\Omega^\perp}$  as  $\Ker P_{\Omega^\perp}=\Im P_{\Omega} $. With some linear algebra, we calculate  $\tau(E_0)=1.2>1$ with $E_0=1/9 X_0$. Moreover, by using the cvx package to solve the spectral norm optimization \eqref{SC22},  the Source Nondegenerate $\rho(E_0)$ is exactly $1$ and gives us a solution $Z_0$. Thus $X_0$ is not a sharp solution of problem \eqref{MC1}. 

However, note further that $Y_0=Z_0+E_0\in \Ker P_{\Omega^\perp}=\Im P_{\Omega} $ is a dual certificate, which is computed by
\[
Y_0=\begin{pmatrix}0&0&1\\0&1&0\\1&0&0\end{pmatrix}. 
\]
Let us check the condition 
\begin{equation}\label{IN1}
\Ker \Phi \cap T_{N_\B(Y_0)}(X_0)=\{0\}
\end{equation}
in \eqref{ma2}. It follows from \eqref{subdif} that 
\[
Y_0=UU^TY_0VV^T=U\begin{pmatrix}1&0&0\\0&0&1\\0&1&0\end{pmatrix}V^T. 
\]
The SVD of the submatrix $\begin{pmatrix}0&1\\1&0\end{pmatrix}$ is certainly \[
\begin{pmatrix}0&1\\1&0\end{pmatrix}\begin{pmatrix}1&0\\0&1\end{pmatrix}\begin{pmatrix}1&0\\0&1\end{pmatrix}.\]
According to \eqref{OY2}, $(\Bar U,\Bar Y)\in \mathcal{O}(X_0)\cap\mathcal{O}(Y_0)$ with 
\[
\Bar U\eqdef\frac{1}{3}\begin{pmatrix}2&1&-2 \\1&2&2\\2&-2&1\end{pmatrix} \quad \mbox{and}\quad \Bar V\eqdef V. 
\]
By formula \eqref{Tang2}, 
\[
T_{N_\B(Y_0)}(X_0)=\left\{\Bar U\begin{pmatrix} A&B\\B^T&C \end{pmatrix}\Bar V^T|\; A\in \R^{1\times 1}, B\in \R^{1\times 2}, C\in \mathbb{S}^2_+\right\}. 
\]
It follows that 
\[
\mathcal{E}=\left\{\Bar U\begin{pmatrix} A&B\\B^T&C \end{pmatrix}\Bar V^T|\; A\in \R^{1\times 1}, B\in \R^{1\times 2}, C\in \R^{2\times 2}\right\}\;\mbox{and}\; \mathcal{E}^\perp=\left\{\Bar U\begin{pmatrix} 0&B\\-B^T&0 \end{pmatrix}\Bar V^T|\;  B\in \R^{1\times 2}\right\}.
\]
To verify \eqref{IN1}, we compute  $\zeta(E_0)$, which is the optimal solution of problem \eqref{Spec2}
\begin{equation*}
    \min_{Z\in \R^{2\times 2}, W\in \R^{3\times 3}}\quad \|Z\|\quad\mbox{subject to}\quad P_{\O^c}(\OU_JZV_J^T+W)=-P_{\O^c}(E_0)\quad \mbox{and}\quad W\in \mathcal{E}^\perp. 
\end{equation*}
This is a convex optimization problem. By using cvx to solve it, we obtain that $\zeta(E_0)=1/6$. As $\zeta(E_0)<1$, $X_0$ is a strong solution of problem \eqref{MC11}.  \endproof
}
\end{Example}

The idea of checking strong solution numerically by using Theorem~\ref{SSN} in the above example will be summarized again in Section~6 with bigger size of nuclear norm minimization problems.

 In the later Corollary~\ref{CoroLR}, we will show a big class of nuclear minimization problem satisfies both Strict Restricted Injectivity \eqref{SS2} and Nondegenerate Source Condition \eqref{NSC2}, but not Restricted Injectivity \eqref{RI}.
 
 Now let us consider a special case of problem \eqref{NNM}
\begin{equation}\label{LRR}
    \min_{X\in \R^{n_1\times n_2}}\quad \|X\|_*\quad\mbox{subject to}\quad LX=M_0,
\end{equation}
where $L$ and $M_0$ are known $q\times n_1$ and $q\times n_2$ matrices, respectively. This is usually referred as the low-rank representation problem \cite{LLYS12}. It is well-known that the optimal solution to problem \eqref{LRR} is unique and determined by $L^\dag M_0$, where $L^\dag$ is the Moore-Penrose opeator of $L$. In the following result, we advance this result by showing that this unique solution is indeed a strong solution, but not necessarily a sharp solution. Indeed, in this certain class,  Strict Restricted Injectivity \eqref{SS2} and Nondegenerate Source Condition \eqref{NSC2} are satisfied, but Restricted Injectivity \eqref{RI} is not.

\begin{Corollary}[Strong minima of low-rank representation problems]\label{CoroLR} Let $L$ be an $q\times n_1$ matrix. If the linear system $LX=M_0$ is consistent, then Strict  Restricted Injectivity \eqref{SS2} and Nondegenerate Source Condition \eqref{NSC2} hold at $X_0\eqdef L^\dag M_0$ in problem \eqref{LRR}.

Consequently, $X_0$ is the strong solution of the low-rank representation problem \eqref{LRR}.    
\end{Corollary}
\noindent{\bf Proof.} Suppose that $U\Sigma V^T$ is a compact SVD to the matrix $L$. Thus $L^\dag=V\Sigma^{-1} U^T$ and  note that  $L^\dag M_0=V\Sigma^{-1} U^TM_0$. Let  $U_0\Sigma_0 V_0^T$ be a compact SVD of $\Sigma^{-1} U^TM_0$ with $\Sigma_0\in \R^{r\times r}$.  Note that 
\[
(VU_0)^TVU_0=U_0^TV^TVU_0=U_0^TU_0=\Id. 
\]
It follows that $VU_0\Sigma V_0^T$ is a compact SVD of $X_0$. By Lemma~\ref{Lem}, we have $E_0\eqdef VU_0 V_0^T\in \ri \partial \|X_0\|_*$. Observe further that 
\[
E_0=VU_0V_0^T=V\Sigma U^T U \Sigma^{-1}U_0V_0^T=L^TU \Sigma^{-1}U_0V_0^T=\Phi^*(U \Sigma^{-1}U_0V_0^T),
\]
 which implies that $E_0\in \Im \Phi^*\cap\ri \partial \|X_0\|_*. $ This verifies Nondegenerate Source Condition \eqref{NSC2} and shows that  $X_0$ is an optimal solution of problem \eqref{LRR}.
 
 Next, let us check Strict  Restricted Injectivity \eqref{SS2}. For any $W\in \Ker \Phi\cap (VU_0)\mathbb{S}^rV_0^T$ with $r=\rank (X_0)$, we find some $A\in \mathbb{S}^r$ such that  $W=VU_0AV_0^T$. We have  
\[
0=\Phi(W)=LW=U\Sigma V^T VU_0AV_0^T=U\Sigma U_0AV_0^T.
\]
It follows that 
\[
0=U_0^T\Sigma^{-1}U^T(U\Sigma U_0AV_0^T)V_0=U_0^T\Sigma^{-1}\Sigma U_0A=U_0^TU_0A=A,
\]
which also implies that $W=0$. This verifies Strict  Restricted Injectivity \eqref{SS2}.

Consequently, according to Corollary~\ref{Ri}, $X_0$ is the strong solution of problem \eqref{LRR}. 
 \endproof
 
 In the following simple example, we show that the unique solution to \eqref{LRR} may be not a sharp solution in the sense \eqref{Sharp}. 
\begin{Example}[Unique solutions of low-rank representation problems are not sharp] {\rm Consider the following optimization problem 
\begin{equation}\label{Sim}
\min_{X\in \R^{2\times 2}}\qquad \|X\|_*\qquad\mbox{subject to}\qquad \begin{pmatrix} 1&1\end{pmatrix}X=\begin{pmatrix}1&0\end{pmatrix}. 
\end{equation}
As $L=\begin{pmatrix}1 &1\end{pmatrix}$ and $M_0=\begin{pmatrix}1&0\end{pmatrix}$, the unique solution to problem \eqref{Sim} is \[
X_0=L^\dag M_0=\begin{pmatrix}0.5\\0.5\end{pmatrix}\begin{pmatrix}1&0\end{pmatrix}=\begin{pmatrix}0.5&0\\0.5&0\end{pmatrix}. 
\]
Pick any $X_\ve\eqdef\begin{pmatrix}0.5+\ve&0\\0.5-\ve&0\end{pmatrix}$ and note that $X_\ve$ satisfies linear constraint in \eqref{Sim}. We have  
\[
\|X_\ve\|_*=\sqrt{(0.5+\ve)^2+(0.5-\ve)^2}=\sqrt{0.5+2\ve^2}. 
\]
It follows that 
\[
\dfrac{\|X_\ve\|_*-\|X_0\|_*}{\|X_\ve-X_0\|_F}=\dfrac{\sqrt{(0.5+\ve)^2+(0.5-\ve)^2}-\sqrt{0.5}}{\sqrt{2\ve^2}}=\dfrac{\sqrt{0.5+2\ve^2}-\sqrt{0.5}}{\sqrt{2\ve^2}}=\dfrac{\sqrt{2}\ve}{\sqrt{0.5+2\ve^2}+\sqrt{0.5}}.
\]
This shows that $X_0$ is not a sharp solution in the sense of \eqref{Sharp}. The hidden reason is that Restricted Injectivity \eqref{RI} is not satisfied in this problem. \endproof
}
\end{Example}

A relatively close problem to \eqref{LRR} is  
\begin{equation}\label{LRRE}
\min_{X\in \R^{n_1\times n_2}}\quad h(LX)+\mu\|X\|_*,
\end{equation}
where the function $h:\R^{q\times n_2}\to [0,\infty]$ satisfies the standing assumptions in Section~4 with open domain and $\mu$ is a positive number. This is a particular case of \eqref{Nucl}. Next we also show that strong minima occurs in this problem.

\begin{Corollary}\label{Redu}
Problem \eqref{LRRE} has a unique and strong  solution.
\end{Corollary}
\noindent{\bf Proof.} It is easy to see that problem \eqref{LRRE} has at least an optimal solution $\OX$.  Let $U\Sigma V^T$ be a compact SVD of $L$ and define   
\[
\OY\eqdef-\dfrac{1}{\mu}L^T\nabla h(L\OX)=-\dfrac{1}{\mu}V\Sigma U^T\nabla h(L\OX)\in\partial \|\OX\|_*. 
\]
Let $U_1\Sigma_1V_1^T$ be a compact SVD of $-\dfrac{1}{\mu}\Sigma U^T\nabla h(L\OX)$ with $\Sigma_1\in \R^{p\times p}$. As $(VU_1)^T(VU_1)=\Id$, it follows that  $VU_1\Sigma_1V_1$ is a compact SVD of $\OY$.  By \eqref{Inver}, we can find $\Bar A\in \mathbb{S}_+^p$ such that $\OX=VU_1\overline{A}V_1^T$, which means 
\[
\overline{A}=U^T_1V^T\OX V_1.
\]
Next let us estimate $T_{N_\B(\OY)}(\OX)$. For any $W\in T_{N_\B(\OY)}(\OX)$, we find sequences $t_k\dn 0$ and $W_k\to W$ satisfying $\OX+t_kW_k\in N_\B(\OY)\subset VU_1\mathbb{S}_+^pV_1^T$ by Lemma~\ref{Lem} again. It follows that \[
W_k\in \frac{1}{t_k}VU_1(\mathbb{S}^p_+-\overline{A}) V_1^T\subset VU_1\mathbb{S}^p V_1^T,
\]
which implies that $T_{N_\B(\OY)}(\OX)\subset VU_1\mathbb{S}^pV_1^T.$ 
 We claim next that $\Ker \Phi \cap T_{N_B(\OY)}(X_0)=\{0\}$. Indeed, take any $W\in T_{N_\B(\OY)}(X_0)$ with $\Phi(W)=0$, we find $B\in \mathbb{S}^p$ such that $W=VU_1BV_1^T$ and 
\[
0=\Phi(W)=LW=U\Sigma V^TVU_1BV_1^T=U\Sigma U_1 B V_1^T.
\]
Hence we have
\[
U_1 B V_1^T=\Sigma^{-1}U^TU\Sigma U_1 B V_1=0.
\]
This implies that $W=0$ and verifies the claim. By Corollary~\ref{CoNu}, $X_0$ is the strong solution of  problem \eqref{LRRE}.\endproof

\section{Numerical Experiments}
\setcounter{equation}{0}

In this section, we perform numerical experiments to demonstrate strong minima, sharp minima, and solution uniqueness for nuclear norm minimization problem \eqref{Nu0}.  The experiments were conducted for different matrix ranks $r$ and numbers of measurements $m$ for $M_0$. Through the section, we also discuss how to use our conditions to check strong minima for problem \eqref{Nu0}.

\subsection{Experiment 1}

In the first experiment, we generate $X_0$, an $n\times n$ matrix of rank $r$, by sampling two factors $W\in\mathbb R^{n\times r}$ and $H\in\mathbb R^{n\times r}$ with independent and identically distributed (i.i.d.) random entries and setting $X_0 = WH^*$. We vectorize problem \eqref{Nu0} in the following form:
\begin{equation}\label{GOP}
    \min_{X\in \R^{n\times n}}\quad \|X\|_*\quad \mbox{subject to}\quad \Phi\ \text{vec}(X) = \Phi\ \text{vec}(X_0),
\end{equation}
where $\Phi\in\mathbb R^{m\times n^2}$ is drawn from the standard Gaussian ensemble, i.e., its entries are independently and identically distributed from a zero-mean unit variance Gaussian distribution. We declare $X_0$ to be recovered (a solution to \eqref{GOP}) if $\|X_{\text{opt}} - X_0\|_F/\|X_0\|_F < 10^{-3}$, as proposed in \cite{CR09}. In order to check sharp minima, it is required to verify Restricted Injectivity \eqref{RI}, compute $\tau(X_0)$ \eqref{tau} or  Source
 Coefficient $\rho(X_0)$ \eqref{SC22}; see \cite{FNT21} or Remark~\ref{Sharp2}. To verify strong minima at $X_0$, we use Strong Source Coefficient $\zeta(X_0)$ in \eqref{Spec} or \eqref{Spec2} and use it . 
 

Specifically, let $U_0\begin{pmatrix}\Sigma_0&0\\0&0\end{pmatrix} V_0^T$ be a full SVD of $X_0$. We respectively denote $u_i$ and $v_j$ as the $i$th and $j$th column of $U_0$ and $V_0$ for $1\le i,j\le n$. Note from the formula of the tangent model space $\TT_0$ in \eqref{T0} that $\mathcal{B}=\{u_iv_j^T: (i,j) \notin [n-r+1,n]\times[n-r+1,n]\}$ forms a basis for $\TT_0$. 
Thus, the Restricted Injectivity holds if $\rank \Phi B = r(2n-r)$, where $B$ is a matrix whose columns are all vectors from the basis $\mathcal{B}$.  

To compute $\tau(X_0)$, we assume that $U\Sigma V^T$ is an SVD of $\Phi$ and denote $V_G$ by the matrix whose columns are the last $n^2-r$ columns of $V$.  We 
 then solve the following vectorized problem for the optimal solution $Z^*$ by using the cvxpy package and compute $\tau(X_0) = \|Z^*\|$:
\begin{equation}\label{TAUOP}
 \min_{Z\in\TT_0^\perp}\|Z\|_F\quad \mbox{subject to}\quad N\text{vec}(Z) = - N\text{vec}(E_0),
\end{equation} 
where $N = V_G^T$ and $\TT_0^\perp$ is known by
\begin{equation}\label{TPERP}
    \TT_0^\perp=\left\{U_0\begin{pmatrix}0&0\\0&D\end{pmatrix}V_0^T|\; D\in \R^{(n-r)\times (n - r)}\right\}.
\end{equation}

To calculate $\rho(X_0)$, we solve the following vectorized problem  of \eqref{SC22} for the optimal value $\rho(X_0)$ by using the cvxpy package:
\begin{equation}\label{RHOP}
 \min_{Z\in\TT_0^\perp}\|Z\|\quad \mbox{subject to}\quad N\text{vec}(Z) = - N\text{vec}(E_0),
\end{equation}
where $N$ and $\TT_0^\perp$ are respectively determined in \eqref{TAUOP} and \eqref{TPERP}. $X_0$ is a sharp solution of problem \eqref{GOP} if either $\tau(X_0)$ or $\rho(X_0)$ is smaller than $1$; see Remark~\ref{Sharp2}. Due to the possible (small) error in computation, we classify sharp minima if $X_0$ is recovered and either $\tau(X_0)<0.99$ or $\rho(X_0)<0.95$.

To classify strong minima, we consider the cases when $X_0$ is recovered, and $\tau(X_0)>0.99$ and $0.95<\rho(X_0)<1.05$. Let $Z_0$ be an optimal solution problem of  \eqref{RHOP} expressed by the following form:
\begin{equation}\label{YZ}
Z_0 = U_0\begin{pmatrix}0&0\\0&D_0\end{pmatrix}V_0^T\quad \mbox{and}\quad Y_0=U_0\begin{pmatrix}\Id&0\\0&D_0\end{pmatrix}V_0^T
\end{equation}
with some $D_0 \in \R^{(n-r)\times n - r}$. Note that $Y_0$ is a dual certificate of $X_0$. According to Theorem~\ref{SSN}, $X_0$ is a strong solution provided that 
\begin{equation}\label{IntP}
\Ker \Phi \cap T_{N_\B(Y_0)}(X_0)=\{0\}. 
\end{equation}
By Theorem~\ref{TheoNu}, this condition holds when the Restricted Injectivity is satisfied and the Strong Source Coefficient $\zeta(X_0)$, the optimal value of problem \eqref{Spec} or \eqref{Spec2} is smaller than $1$. Let  $\Hat U\Hat \Sigma\Hat V$ be an SVD of $D_0$. We write $U_0=[U_I\; U_J]$ and $V_0=[V_I\;V_K]$, where $U_I$ and $V_I$ are the first $r$ columns of $U_0$ and $V_0$, respectively. Define $\OU = [U_I\; U_J\Hat U]$ and $\OV = [V_I\; V_J\Hat V]$, it follows that  $(\OU,\OV)\in \mathcal{O}(X_0)\cap\mathcal{O}(Y_0)$ by \eqref{OY2}.   To compute $\zeta(X_0)$, we solve the vectorized problem of \eqref{Spec2}:
\begin{equation}\label{ZETAP}
 \min_{Z\in\TT_0^\perp, W\in\mathcal E^\perp}\|Z\|\quad \mbox{subject to}\quad N\text{vec}(Z + E_0 + W) = 0,
\end{equation}
where $\TT_0^\perp$ is determined in \eqref{TPERP}, $E_0=U_IV_I^T$, and $\mathcal E^\perp$ is taken from \eqref{Eperp}:
\begin{equation}\label{EPERP}
    \mathcal{E}^\perp=\left\{\OU\begin{pmatrix}A&B&C\\-B^T&0&0\\D&0&0\end{pmatrix}\OV^T|\;A\in \mathbb{V}_r, B\in \R^{r\times (p-r)}, C\in \R^{r\times (n-p)}, D\in \R^{ (n-p)\times r} \right\}.
\end{equation}
 We classify strong (non-sharp) minima if  $\zeta(X_0) < 0.95$. 
 
\begin{figure*}[!htpb]
\vspace{-1ex}
\begin{center}
\includegraphics[width=1\linewidth]{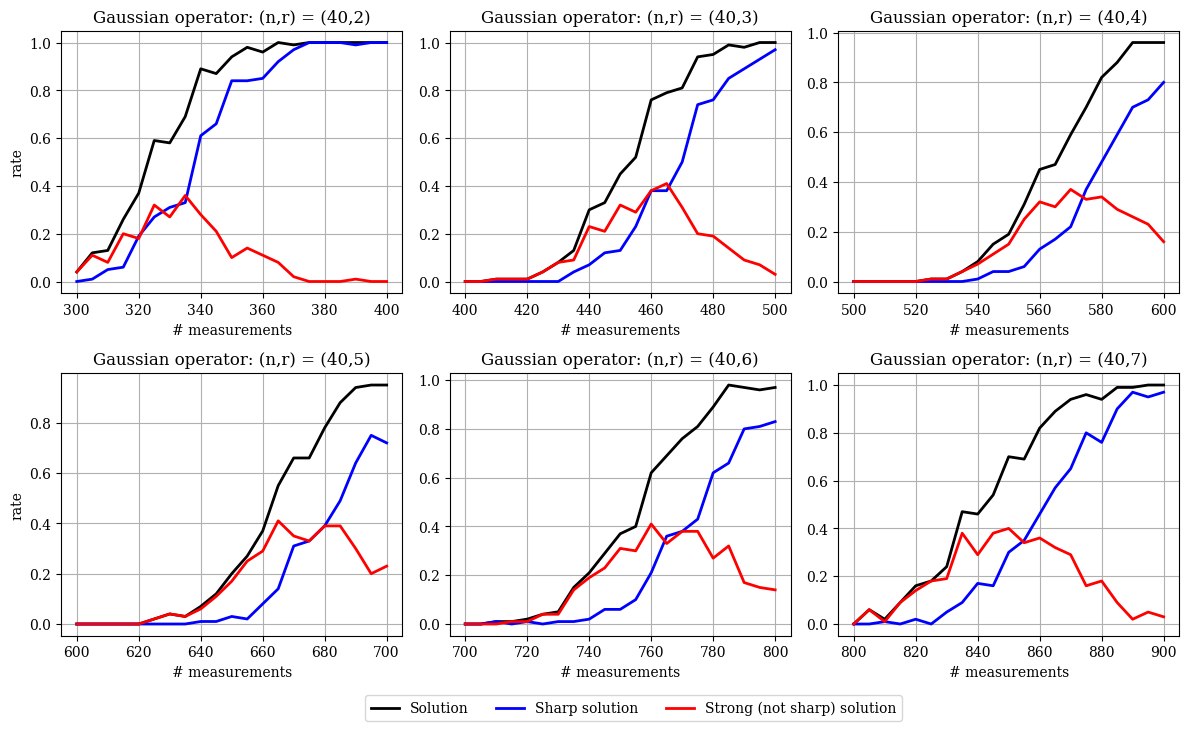}
\caption{Proportions of cases for which $X_0$ is a solution, sharp solution, and strong (not sharp) solution with respect to the number of measurements.\label{fig:GOP}} 
\end{center} 
\vspace{-3ex}
\end{figure*}

 To illustrate the occurrence of strong minima, sharp minima, and solution uniqueness of problem \eqref{GOP}, we create following graphs representing the proportion of each situation with respect to the number of measurements in Figures \ref{fig:GOP}. For each fixed $n=40$ and $2\le r\le 7$, at any measurements $m$  we study 100 random cases and record the percentage of cases that are recovered, sharply recovered, and strongly (not sharply) recovered in black, blue, and red curves, respectively.  Observe that the percentage of cases where $X_0$ is a strong (not sharp) solution is highest at approximately 40\% and is more than the cases of sharp minima when the number of measurements is not big enough. This phenomenon was obtained at different measurements for different ranks, indicating a significant number of cases with strong (not sharp) solutions. Additionally, higher ranks require more measurements to achieve the highest percentage of cases with strong (not sharp) solutions.

\begin{figure*}[!htpb]
\vspace{-1ex}
\begin{center}
\includegraphics[width=1\linewidth]{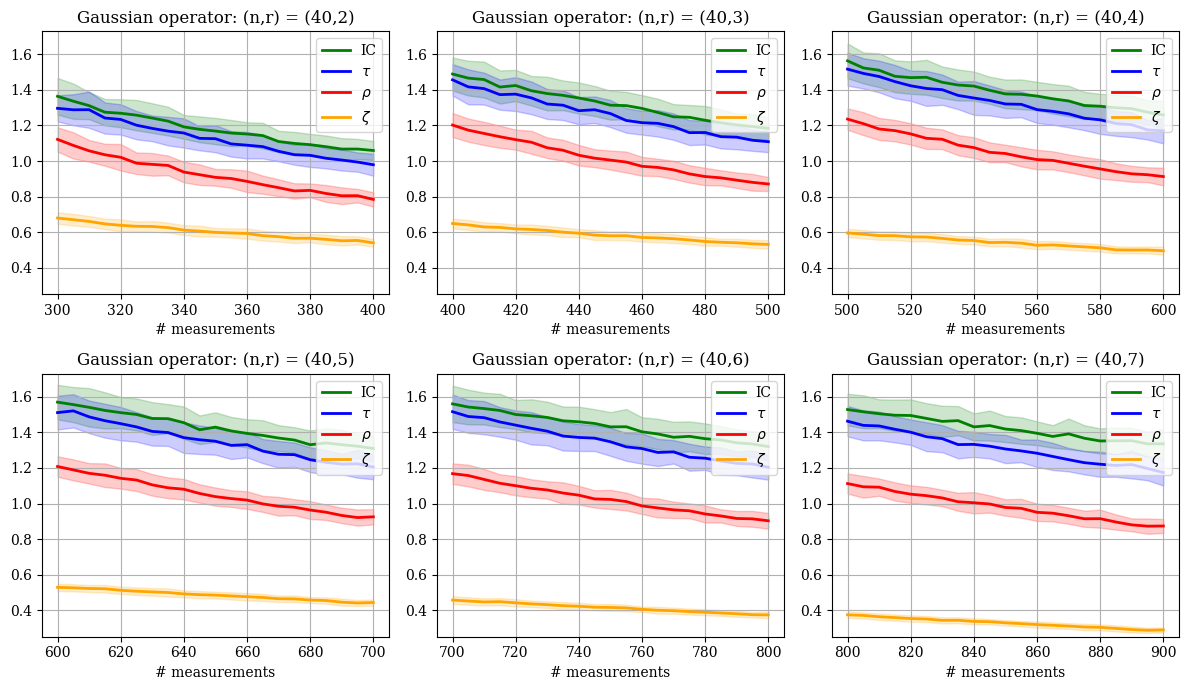}
\caption{Evolution of the average value of $\tau(E_0)$, Source
Nondegenerate $\rho(E_0)$, and Strong Source Coefficient $\zeta(E_0)$ with respect to the number of measurements. \label{fig:GOPvalue}} 
\end{center} 
\vspace{-3ex}
\end{figure*}

 We also plot the average values of $\tau(X_0)$, IC($X_0$) (Irrepresentability Criterion \eqref{IC}), $\rho(X_0)$, and $\zeta(X_0)$ for each number of measurements in Figures \ref{fig:GOPvalue} with different ranks. It seems that using $\rho(X_0)$ to check sharp minima gives us more cases than using $\tau(X_0)$. Moreover, $\zeta(X_0)$ is significantly smaller than both $\tau(X_0)$ and $\rho(X_0)$ while IC($X_0$) is slightly greater than $\tau(X_0)$.   


\subsection{Experiment 2}
In the second experiment, we study the following matrix completion problem:
\begin{equation}\label{MCP}
    \min_{X\in \R^{n\times n}}\quad \|X\|_*\quad \mbox{subject to}\quad X_{ij}=(X_0)_{ij},\quad (i,j)\in\Omega.
\end{equation}
by a similar process to the first one. We also generate $X_0$, an $n\times n$ matrix of rank $r$, by sampling two factors $W\in\mathbb R^{n\times r}$ and $H\in\mathbb R^{n\times r}$ with i.i.d. random entries and setting $X_0 = WH^*$. However, this time we sample an indexed subset $\Omega$ of $m$ entries uniformly at random from $[n]\times[n]$. Cvxpy package is also used to solve problem \eqref{MCP} with an optimal solution $X_{\text{opt}}$.  
 $X_0$ is said to be recovered if $\|X_{\text{opt}} - X_0\|_F/\|X_0\|_F < 10^{-3}$. To classify sharp minima, we check Restricted Injectivity \eqref{RI} (as done in the first experiment), compute $\tau(X_0)$ in \eqref{tau} or  Source
Coefficient $\rho(X_0)$ in \eqref{SC22}, and restrict $\tau(X_0)\le0.99$ or  $\rho(X_0)\le 0.95$.


Specifically, to calculate $\tau(E_0)$, we follow the following steps. It is similar to the first experiment, denote  $u_i$ and $v_j$ by the $i$th and $j$th column of $U_0$ and $V_0$, where $U_0DV_0^T$ is a full SVD of $X_0$. We define $B_{ij} = P_\Omega(u_iv_j^T)$ for all $(i,j)\in\Gamma\eqdef\{(i,j)\in[n]\times[n]|\; (i,j)\notin [n-r,n]\times[n-r,n]\}$, where $P_\Omega$ is the projection mapping defined by $P_\Omega(X)_{ij} = X_{ij}$ if $(i,j)\in\Omega$ and 0 otherwise. Then we solve the following linear system for $\al_{ij}$: 
\[
\sum_{(i,j)\in\Gamma}\alpha_{ij}P_\Gamma\left(U_0^TB_{ij}V_0\right) = \begin{pmatrix}I_r&0\\0&0\end{pmatrix}.
\]
It is not difficult to obtain from \eqref{tau} that
\[
\tau(E_0)=\left\|Y - E_0\right\|,
\]
where $Y = \sum_{(i,j)\in\Gamma}\alpha_{ij}B_{ij}$ and $E_0 = U_0\begin{pmatrix}I_r&0\\0&0\end{pmatrix}V_0^T$.

To compute $\rho(X_0)$, we transform problem \eqref{SC22} to the case of matrix completion as follows:  
\begin{equation}\label{RHOMCP}
 \min_{Z\in\TT_0^\perp}\|Z\|\quad \mbox{subject to}\quad Z_{ij} = - (E_0)_{ij},\quad (i,j)\notin \Omega
\end{equation}
and solve by using cvxpy package for the optimal value $\rho(X_0)$,  
where $\TT_0^\perp$ is determined in \eqref{TPERP}.

\begin{figure*}[!htpb]
\vspace{-1ex}
\begin{center}
\includegraphics[width=1\linewidth]{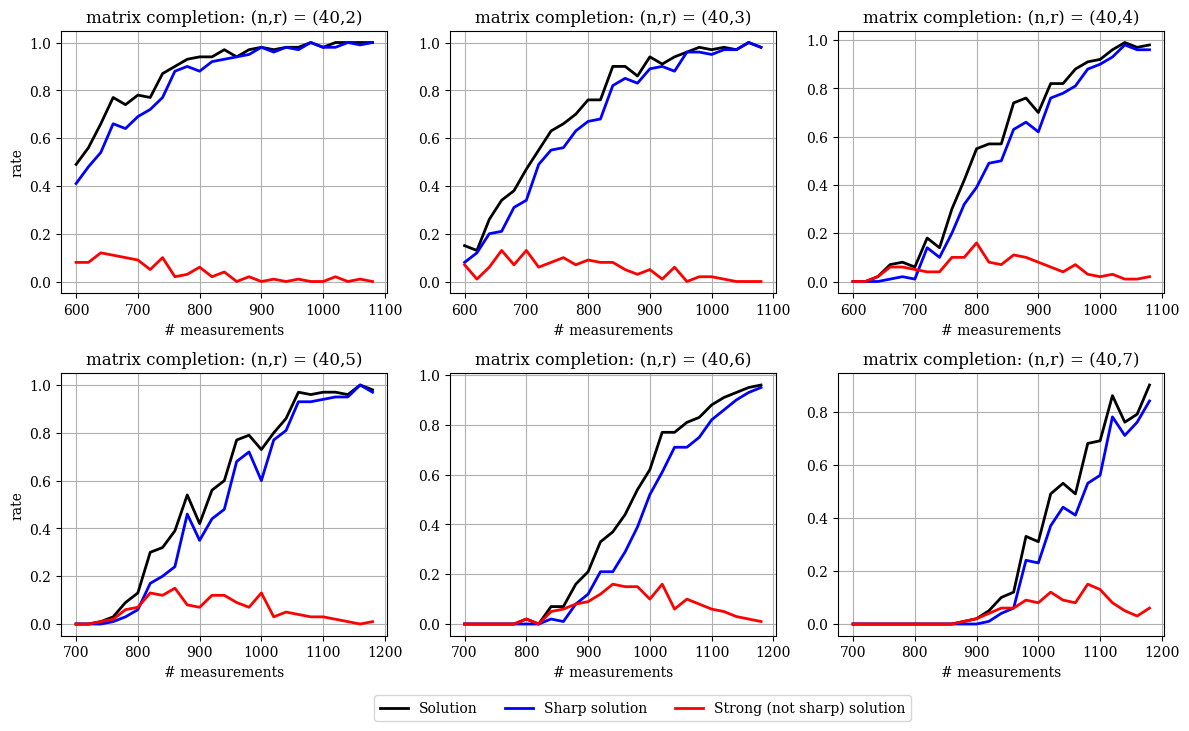}
\caption{Proportions of cases for which $X_0$ is a solution, sharp solution, and strong (not sharp) solution with respect to the number of measurements.\label{fig:MCP}} 
\end{center} 
\vspace{-3ex}
\end{figure*}

\begin{figure*}[!htpb]
\vspace{-1ex}
\begin{center}
\includegraphics[width=1\linewidth]{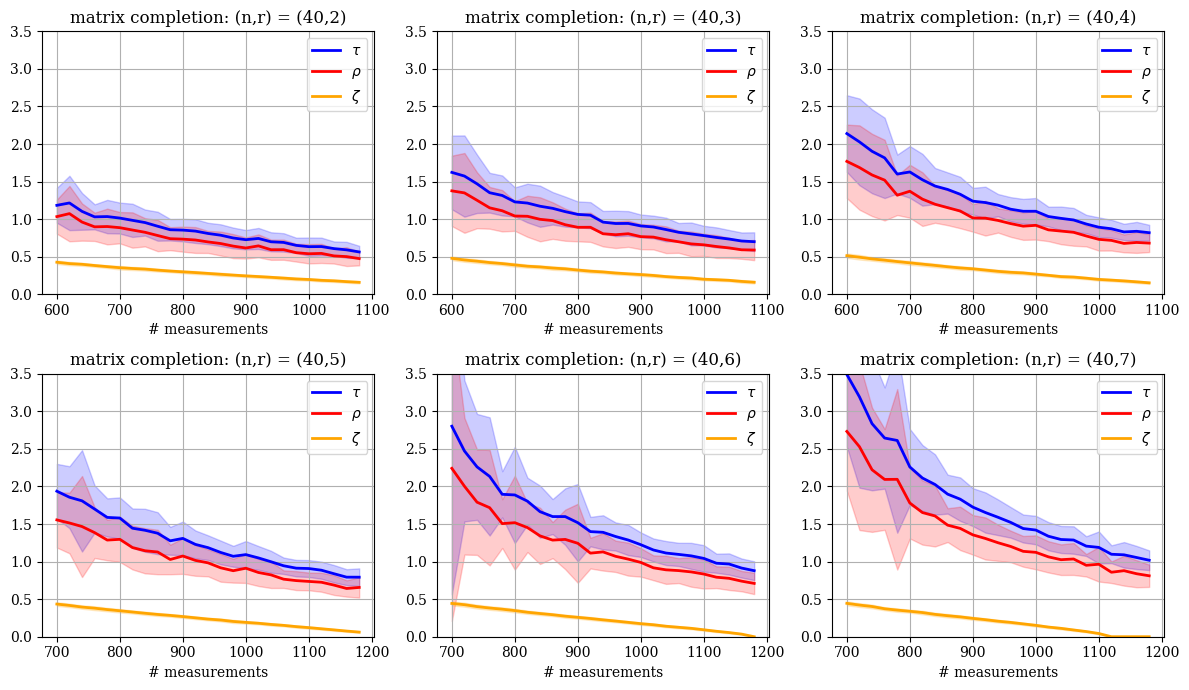}
\caption{Evolution of the average value of $\tau(E_0)$, Source
Nondegenerate $\rho(E_0)$, and Strong Source Coefficient $\zeta(E_0)$ with respect to the number of measurements. \label{fig:MCPvalue}} 
\end{center} 
\vspace{-3ex}
\end{figure*}

Similarly to  \eqref{YZ}, we denote $Z_0$ be an optimal solution of problem \eqref{RHOMCP} and denote $Y_0=Z_0+E_0$ as a dual certificate of $X_0$. To check if $X_0$ is a strong solution, we only need to verify \eqref{IntP} with $\Ker \Phi=\Ker P_\Omega$ by Theorem~\ref{SSN}. According to Theorem~\ref{TheoNu}, the later holds under the Restricted Injectivity and $\zeta(X_0)<1$, which is the optimal solution of the following problem, a version of \eqref{Spec2} for matrix completion: 
\begin{equation}\label{ZETAMCP}
    \min_{Z\in\TT_0^\perp, W\in\mathcal E^\perp}\|Z\|\quad \mbox{subject to}\quad (Z + E_0 + W)_{ij} = 0,\quad (i,j)\notin \Omega.
\end{equation}
Here  $\TT_0^\perp$  and $\mathcal E^\perp$ are defined \eqref{TPERP} and \eqref{EPERP}. We classify a case of strong minima (strong recovery) if $X_0$ is recovered and $\zeta(X_0)\le 0.95$, but $\tau(X_0)>0.99$ and $0.95<\rho(X_0)<1.05$. In Figures \ref{fig:MCP} , we plot the proportion of cases when $X_0$ is a unique solution, sharp solution, and strong (not sharp) solution in relation to the number of measurements. Additionally, Figures \ref{fig:MCPvalue}  displays the curves of the average value of $\tau(X_0)$, $\rho(X_0)$, and $\zeta(X_0)$ with respect to the number of measurements.

Based on Figure \ref{fig:MCP}, we can see that the highest percentage of cases where $X_0$ is a strong solution (but not a sharp one) is just around 15\%. It is much smaller than the 40\% in Experiment 1. A possible reason is due to the special structure of the linear operator $\Phi=P_\Omega$ in \eqref{MCP},  which is chosen in such a way that each row contains only one entry of 1 and the remaining entries are 0. Additionally, these entries of 1 must be in distinct columns. Similar to Experiment 1, we found that higher ranks require more measurements to achieve the highest percentage of cases with strong (not sharp) solutions. 

As shown in Figure \ref{fig:MCPvalue}, the average values of $\tau(X_0)$, $\rho(X_0)$, and $\zeta(X_0)$ change depending on the number of measurements. The difference between the curves of $\tau(X_0)$ and $\rho(X_0)$ is less noticeable in comparison to Experiment 1, but the curve $\zeta(X_0)$ is still significantly lower than both of $\tau(X_0)$ and $\rho(X_0)$. 

\begin{thebibliography}{99}


\bibitem{ALMT14} D. Amelunxen, M. Lotz, M. B. McCoy, and J. A. Tropp: Living on the edge: A geometric theory of phase transitions in convex optimization, \textit{IMA Inform. Inference}, {\bf 3} (2008), 224--294.

\bibitem{AG08} F. J. Arag\'on Artacho and M. H.  Geoffroy: Characterizations of metric regularity of subdifferentials, \textit{J. Convex Anal.}, {\bf 15} (2008), 365--380.

\bibitem{AV11} B. P. W.  Ames and S. A. Vavasis: Nuclear norm minimization for the planted clique and biclique problems, \textit{Math. Program.}, \textbf{29} (2011), 69--89. 

\bibitem{BT80}
A.~Ben-Tal.
\newblock Second order and related extremality conditions in nonlinear programming.
\newblock{\em J. Optim. Theory Appl.}, \textbf{31} (1980) 143--165.

\bibitem{BZ82} A.~Ben-Tal and J. Zowe: A unified theory of first- and second order conditions for extremum problems in topological vector spaces. \textit{Math. Program.} 19 (1982), 39--76.


\bibitem{BI95a} J. F. Bonnans and  A. D. Ioffe: second order sufficiency and quadratic growth for non-isolated minima, \textit{Math. Oper. Res.}, \textbf{20} (1995), 801--817.


\bibitem{BCS99} J.F. Bonnans, R. Cominetti, and A. Shapiro: Second order optimality conditions based on parabolic second order tangent, \textit{SIAM J. Optim.}, \textbf{9} (1999), 466-492

\bibitem{BS00}  J. F. Bonnans and A. Shapiro: {\em Perturbation  Analysis of Optimization Problems}, Springer, 2000.

 \bibitem{BNPS17} J. Bolte, T. P. Nguyen, J. Peypouquet, B. W. Suter: From error bounds to the complexity of first-order descent methods for convex functions, \textit{Math. Program.}, {\bf 165} (2017), 471--507. 
 
 \bibitem{B87} J. Burke: Second order necessary and sufficient conditions for convex composite NDO, \textit{Math. Program.}, {\bf 38} (1987), 287-302.
 
 
 \bibitem{C12} P. M. Cohn: \textit{Universal algebra}, Springer Science \& Business Media, 2012.
 
 
\bibitem{C77} L. Cromme: Strong uniqueness. A far-reaching criterion for the convergence analysis of
iterative procedures, \textit{Numer. Math.}, 29(2) (1977/78), 179--193


\bibitem{CRPW12} V. Chandrasekaran, B. Recht, P. A.  Parrilo, and A. S. Willsky: The convex geometry of linear inverse problems, \textit{Found Comput Math}, \textbf{12} (2012), 805--849.


\bibitem{CP10} E.  Cand\`es
and Y. Plan: Matrix completion with noise, \textit{Proceeding of the IEEE}, \textbf{98} (2010), 925--936.

\bibitem{CR09} E. Cand\`es and B. Recht: Exact matrix completion via convex optimization \textit{Found. Comput. Math.}, \textbf{9} (2009), 717--772.

\bibitem{CR13}  E. Cand\`es and B. Recht: Simple bounds for recovering low-complexity models, \textit{Math. Program.},  \textbf{141} (2013), 577--589.






\bibitem{CD19} Y. Cui and C. Ding: Nonsmooth composite matrix optimization: strong regularity, constraint nondegeneracy and beyond,  \textit{ArXiv},	arXiv:1907.13253

\bibitem{CDZ17} Y. Cui, C. Ding, and X.  Zhao: Quadratic growth conditions for convex matrix optimization problems associated with spectral functions, \textit{SIAM J. Optim.}, {\bf 27} (2017), 2332--2355.

\bibitem{CST19}  Y. Cui, D. Sun,  and K-C. Toh: On the R-superlinear convergence of the KKT residuals generated by the augmented Lagrangian method for convex composite conic programming, \textit{Math. Program.}, {\bf 178} (2019),381--415


\bibitem{D17} C. Ding: Variational analysis of the Ky Fan k-norm, \textit{Set-Valued and Var. Anal.}, \textbf{25} (2017),  265--296. 


\bibitem{DMN14} D. Drusvyatskiy, B. S. Mordukhovich, and T. T. A.  Nghia: second order growth, tilt stability, and metric regularity of
the subdifferential,  \textit{J. Convex Anal.}, {\bf 21} (2014), 1165--1192.

\bibitem{FM68} A. V. Fiacco and G. P. McCormick: \textit{Nonlinear Programming: Sequential Unconstrained
Minimization Techniques}, Wiley, New York, 1968.

\bibitem{FMP18}  J. Fadili, J. Malick, G. Peyr\'e: Sensitivity Analysis for Mirror-Stratifiable Convex Functions, \textit{SIAM J. Optim.}, \textbf{28} (2018), 2975---3000. 

\bibitem{FNT21} J. Fadili, T. T. A. Nghia, and T. T. T. Tran: Sharp, strong and unique minimizers for low complexity robust recovery, \textit{IMA Inf. Inference}, \textbf{12} (2023).



\bibitem{GHS11}  M. Grasmair, M.  Haltmeier, and  O. Scherzer: Necessary and sufficient conditions for linear convergence of $\ell_1$ regularization,  \textit{Comm.  Pure Applied Math.} \textbf{64}(2011), 161--182.


\bibitem{H92} J. Harris: {\em Algebraic Geometry}, Springer, 1992. 

\bibitem{HO14} C. J. Hsieh and P. Olsen: Nuclear norm minimization via active subspace selection. \textit{ICML}, \textbf{32} (2014), 575--583. 

\bibitem{HP23} T. Hoheisel and E. Paquette: Uniqueness in nuclear norm minimization: Flatness of the nuclear norm sphere and simultaneous polarization, \textit{J. Optim. Theo. Appl.}, {\bf 197} (2023), 252-276.

\bibitem{LFP17} J. Liang,  J. Fadili, G. Peyr\'e: Activity identification and local linear convergence of forward-backward-type methods, \textit{SIAM  J. Optim.}, {\bf 27} (2017), 408--437. 

\bibitem{LS05a} A. S. Lewis and H. S. Sendov: Nonsmooth analysis of singular values. I. Theory.
\textit{Set-Valued Anal.}, {\bf 13} (2005), 213--241.

\bibitem{LS05b} A. S. Lewis and H. S. Sendov: Nonsmooth analysis of singular values. II. Applications, \textit{Set-Valued Anal.}, {\bf 13} (2005), 243--264.

\bibitem{LT93} Z.-Q. Luo and P. Tseng: Error bounds and convergence analysis of feasible descent methods: a general approach, \textit{Ann. Oper. Res.}, {\bf 46} (1993), 157--178.

\bibitem{LLYS12}G. Liu, Z. Lin, S. Yan, J. Sun, Y. Yu, and Y. Ma: Robust recovery of subspace structures by low-rank representation. \textit{IEEE Transactions on Pattern Analysis and Machine Intelligence}, \textbf{35} (2012), 171-184.

\bibitem{I79}
A. D. Ioffe.
\newblock Necessary and sufficient conditions for a local minimum III: Second order conditions and augmented duality.
\newblock {\em SIAM J. Control Optim.}, \textbf{17} (1979), 266--286.

\bibitem{I91} A. D. Ioffe: Variational analysis of a composite function: A formula for the lower second
order epi-derivative, \textit{J. Math. Anal. Appl.}, \textbf{160} (1991), 379--405.



\bibitem{M1} B. S. Mordukhovich: {\em Variational Analysis and Generalized Differentiation}, Springer, 2006.


\bibitem{MS23} A. Mohammadi and E. Sarabi: Parabolic regularity of spectral functions. Part I: Theory. 	\textit{arXiv:2301.04240} 



\bibitem{NDEG13} S. Nam, M. E. Davies, M. Elad, R.  Gribonval: The cosparse analysis model and algorithms, \textit{Applied and Computational Harmonic Analysis}, \textbf{34} (2013), 30--56.

\bibitem{P79} B. T. Polyak. Sharp minima. Technical report, Institute of Control Sciences Lecture Notes,
Moscow, USSR, 1979.

\bibitem{RFP10} B. Recht, M. Fazel, and P. Parrilo: Guaranteed minimum rank solutions of matrix equations via nuclear
norm minimization. \textit{SIAM Rev.}, \textbf{52} (2010), 471--501.


 \bibitem{R70} R. T. Rockafellar: {\em Convex Analysis}, Princeton University Press, Princeton, New Jersey, 1970.
 
 \bibitem{R88} R. T. Rockafellar: second order optimality conditions in nonlinear programming obtained by way of epi-derivatives, \textit{Trans. Amer. Soc.}, {\bf 307} (1988), 75--108. 


\bibitem{RW98} R. T. Rockafellar and  R. J-B. Wets: {\em Variational Analysis}, Springer, Berlin, 1998.




\bibitem{SW99} M. Studniarski and D. E. Ward: Weak sharp minima and suffficient conditions, \textit{SIAM J. Control Optim.}, {\bf 38} (1999), 219--236.


\bibitem{VPF17} S. Vaiter, G. Peyr\'e, and J. Fadili: Model consistency of partly smooth regularizers. \textit{IEEE Transactions on Information Theory}, \textbf{64} (2017), 1725--1737.



\bibitem{VGFP15} S. Vaiter, M. Golbabaee, J. Fadili, and G. Peyr\'e: Model selection with low complexity priors. \textit{Information and Inference: A Journal of the IMA}, \textbf{4}  (2015), 230--287.

\bibitem{VPF15} S. Vaiter, G. Peyr\'e, and J. Fadili: Low complexity regularization of linear inverse problems. In {\em Sampling Theory, a Renaissance}, 103--153. Birkh\"auser, Cham, 2015.


\bibitem{W92} G. A. Watson: Characterization of the subdifferential of some matrix norms, \textit{Linear Algebra Appl.}, {\bf 170} (1992), 33--45

\bibitem{W94}
D. E. Ward: Characterizations of strict local minima and necessary conditions for
  weak sharp minima, {\em J. Optim. Theory and Appl.},
  \textbf{80} (1994), 551--571.

\bibitem{ZS17} Z. Zhou and A. M-C.  So: A unified approach to error bounds for structured convex optimization, \textit{Math. Program.}, {\bf 165} (2017), 
689--728.


\bibitem{ZZX13} L. Zhang, N. Zhang, and X.  Xiao: On the second order directional derivatives of singular values of matrices and symmetric matrix-valued functions, \textit{Set-Valued  Var. Anal.}, {\bf 21} (2013), 557--586.

\bibitem{ZT95}
R. Zhang and J. Treiman:
 Upper-Lipschitz multifunctions and inverse subdifferentials, \textit{\em Nonlinear Anal.}, {\bf 24} (1995), 273--286.
 

\bibitem{WYYZZ22}
X.~Wang, J.~J. Ye, X.~Yuan, S.~Zeng, and J.~Zhang.: Perturbation techniques for convergence analysis of proximal gradient
  method and other first-order algorithms via variational analysis, \textit{Set-Valued and Variational Analysis}, {\bf 30} (2022), 39--79.

\end{thebibliography}
\end{document}